\documentclass[a4paper, 12pt]{article}

\usepackage{amsmath,amssymb,amsfonts,setspace}
\usepackage{mathtools}
\usepackage{graphicx}
\usepackage{authblk}
\usepackage{pgfplots}\pgfplotsset{compat=1.17}

\newcommand{\eq}{\begin{equation}}
\newcommand{\en}{\end{equation}}
\newcommand{\prob}{{\mathbb P}}
\newcommand{\me}{{\mathbb E}}

\newtheorem{theorem}{Theorem}
\newtheorem{lemma}{Lemma}

\newtheorem{rem}{Remark}
\newtheorem{prop}{Proposition}

\def\endpf{\hfill $\Box$ \vskip0.5cm}
\def \proof{\noindent{\it Proof.\ }}

\newlength{\bibitemsep}
\newlength{\bibparskip}
\let\oldthebibliography\thebibliography
\renewcommand\thebibliography[1]{%
  \oldthebibliography{#1}%
  \setlength{\parskip}{\bibitemsep}%
  \setlength{\itemsep}{\bibparskip}%
}

\begin{document}

\title{Optimal Stopping for the Uniform Distribution}
\author{Alexander  Gnedin}
\maketitle

\begin{abstract}
\noindent
Many discrete-time optimal stopping problems are known to have more tractable limit forms based on a planar Poisson process.
Using this tool we find a solution to the optimal stopping problem for  i.i.d. sequence
of $n$ discrete  uniform random variables, in  the asymptotic regime where $n$  and the range of distribution are of the same order.
The optimal stopping rule in the Poisson problem is identified, by means of a time change,   with  known asymptotic solution  to Lindley's problem of minimising the expected rank.
\end{abstract}

\section{Introduction}
The following problem  introduced by Moser \cite{Moser} belongs to classics of  the discrete-time optimal stopping.
A  number is drawn at random  from the uniform distribution on $[0,1]$. We may either keep the number, or
reject it and draw again. We can then either keep the second
number  drawn or reject it, and draw again, and so on until being willing to accept the last  observed number.
Suppose we have at most $n$ opportunities of which exactly one must be used to eventually accept a number.
Let $V_n$ be the smallest  possible expected value of the accepted draw. What is $V_n$ and  what is a stopping rule that achieves this value?

Contrary to the established  tradition we formulated  the problem
as a {\it minimisation} task, to cast the principal  large-$n$ asymptotics
in the most  transparent way
as
\begin{equation}\label{Mos}
\lim_{n\to\infty} nV_n=2,
\end{equation}
and to facilitate linking with the scaling limits in the sequel.
Extensions of Moser's problem to sampling from other  continuous distributions (possibly  changing  from trial to trial) have been widely studied, see \cite{Katriel} for recent developments and bibliography.
The related vast area of prophet inequalities is focused on comparison of the optimal stopping value with
expectation of the most favourable draw value that would be
 accessible under a complete foresight of  random sequence \cite{Malin}.

Historically, the first problem of optimal stopping
was proposed by Cayley \cite{Cayley} for sampling without replacement from an urn of lottery tickets, and solved explicitly by backward induction for
payoffs $\{1,2,3,4\}$ and
$n=1,2,3,4$ draws.
Ferguson in his  textbook   \cite{Ferguson} regards Moser's problem as approximation to  Cayley's model for sampling from $\{1,2,\ldots,N\}$ when $N$ is large. However,  if both $N$ and $n$ are large
the approximation can only be justified  when $N$ is  of the order  higher than   $n$,  to ensure that
  each  particular value appears among $n$ draws with negligible probability.
If   $N$ and $n$ are  of comparable magnitude
one should be careful
to distinguish between sampling with and without replacement. Sampling without replacement is  always more favourable
\cite{Brien} but  gives rise to a much more complex stopping problem apparently unsolved up to date (for $n<N$ of the same order).

In this paper we first revise Moser's problem with the accent on its limit form based on the planar Poisson process,  adding some insights to the existing literature \cite{GnedinSPA,  KK1, Rusch}.
Then we move into the unexplored territory to consider the optimal stopping for sampling with replacement from the discrete uniform distribution as  $n/N\to T, 0<T<\infty$.
This asymptotic regime  is known as the central domain of parameters in the classic  balls-in-boxes  allocation scheme  \cite{Kolchin}.
We find the limit stopping value function and the optimal strategy, in particular show that for $n=N$ the constant analogous to  (\ref{Mos}) is about $2.513$.
The convergence of the discrete stopping values to the Poisson limit is proved by an elementary argument, which interprets the optimality recursion as an Euler scheme for the limit differential equation.
 Finally we find out, quite unexpectedly,  that the optimal strategy transforms, via a simple time change, into known  asymptotic solution to
 Lindley's  problem of   minimising the expected  rank \cite{CMRS, Ferguson, Lindley, Mucci2}; at the level of the Poisson limits the connection is an equivalence of  two stopping problems.
The analytical aspects of the exact solution and two-sided  bounds enrich  the bunch of examples of  stopping boundaries for Poisson processes found in \cite{Rusch}.

Though we introduce all necessary tools, skipping some routine details  will be inevitable  to keep the exposition concise.
See \cite{Resnick} for the general background on point processes and extreme values, and
 \cite{CRS, Ferguson, Peskir} for basics of the optimal stopping theory.

\section{Moser's problem and its limit form}

\subsection{Moser's problem in a nutshell}
Suppose the successive draw values     are  independent random variables $U_1, U_2,\ldots$ with $[0,1]$-uniform distribution. We define the stopping value for the problem with (at most) $n$ trials as
the `minimal expected loss'
\begin{equation}
\label{ContiMP}
V_n:=\inf_{\sigma\leq n} {\mathbb E}\, U_\sigma,
\end{equation}
where the infimum is taken over all stopping times $\sigma$ adapted to the sequence and taking values in $\{1,\ldots,n\}$.
In this generality, at  each step the observer may use the information contained in all draws available so far.
By the independence assumption,  however,
in the event  that stopping has not occurred before step $j<n$, only $U_j$ and the number of  remaining trials  $n-j$ do matter for deciding between accepting  the $j$th draw hence terminating, or continuing with the next draw.

The optimality principle  dictates to resolve the stopping-continuation dilemma at each step by comparing the immediate loss of stopping with the best  expected outcome by continuation
(given the observed data).
For the problem (\ref{ContiMP})
this implies the  Wald-Bellman optimality equation  $V_n=\me \min(U_1, V_{n-1})$, $n\geq 1$, which can be re-arranged as
$$V_n-V_{n-1}=-\me \,(V_{n-1}-U_1)_+$$
and upon integration   becomes
\begin{equation}\label{BE}
V_n-V_{n-1}=-\frac{1}{2}\, V_{n-1}^2, ~~~~n\geq 1,~~~V_0:=1.
\end{equation}
Relative to the temporal succession of draws, the  difference equation (\ref{BE}) is a {\it backward} recursion for the {\it continuation value} $V_n$,
which assesses the best outcome for an online  decision maker
when $n$ steps remain to go.

Accordingly, the optimal stopping time is of the threshold form
\begin{equation}\label{n-sigma}
\sigma_n:=\min\{1\leq j\leq  n: U_j\leq V_{n-j}\},
\end{equation}
where the initial value $V_0=1$ forces acceptance of the last draw (since $U_n\leq 1$), while (\ref{BE}) with $n=1$ gives $V_1=\me\, U_1=1/2$.
It is  a peculiar feature of the uniform distribution that the same number $V_{n-j}$ carries  the twofold role
of a critical threshold for draw $U_j$ to be accepted and  the optimal continuation outcome when   $n-j$ trials  are still  available.

The limit (\ref{Mos}) follows from  the more informative bounds on the stopping value, which  adjusted to the minimisation problem (\ref{ContiMP}) are
\begin{equation}\label{MoserBounds}
\frac{2}{n+H_n+2}\leq  V_n < \frac{2}{n+H_n+1},
\end{equation}
where $H_n:=\sum_{k=1}^n 1/k$ is the harmonic number.  See
 \cite{Ferguson, Finch, Moser}  for this and more precise asymptotic approximations.
The stopping value  can be compared
with the   best possible  value of a prophet
$$n\, \me \min (U_1,\ldots,U_n)=\frac{n}{n+1}\to 1.$$

We note in passing that (\ref{BE}) is the recursion for the survival probabilities of the critical Galton--Watson process with binary splitting, in which every particle has $0$ or $2$ offspring with probability $1/2$ each.

The following asymptotic moments of the optimal stopping time
were  obtained  in \cite{Ent, Mazalov} by  involved combinatorial  calculations based on   (\ref{n-sigma}):
\begin{equation}\label{MM}
\frac{1}{n}\,\me\,\sigma_n\to\frac{1}{3},~~~\frac{1}{n^2}\,{\rm Var}(\sigma_n)\to \frac{1}{18} ~~~~~({\rm as}~~n\to\infty).
\end{equation}
Matching the moments  with  (\ref{MM}), one might expect that the distribution of  $\sigma_n/n$ is close to  beta$(1,2)$-distribution on $[0,1]$.
Confirming the guess and deriving the limit law for the accepted draw $U_{\sigma_n}$ are most naturally done in the continuous-time scenario to follow.

%
%

\subsection{The planar Poisson framework}\label{S2.2}

The adjustment of Moser's problem  to sampling from the uniform distribution on $[0,T]$ or any other interval is obvious. Nevertheless
 introducing a scaling parameter  will be beneficial to reveal a trade-off between the  range of draws and horizon of the stopping problem.

The random scatter of points on the plane
\begin{equation}\label{SP}
\left\{\big((T/n)j,\,(n/T)U_j\big), ~j\in {\mathbb N} \right\}
\end{equation}
for large $n$ can be approximated  by  a Poisson point process $\Pi$  with unit rate
 in the positive quadrant  ${\mathbb R}_+^2$. Simply put,  the number of points  (\ref{SP}) falling in a bounded domain $D\subset {\mathbb R}_+^2$ has approximately Poisson distribution with expectation
equal to the area of $D$,
and  the point counts over nonoverlapping domains are almost independent.
Formally, the approximation is captured  by the measure-theoretic   concept of `vague convergence'   of point random measures \cite{Resnick}, that
for (\ref{SP})
 takes care of $U_j$'s of the order $n^{-1}$ and times $j$ of the order $n$,  while wiping out to infinity all significantly larger or later draws.

For the first $n$ uniform draws, with indices and values scaled as  in (\ref{SP}), we will have  the  bivariate point scatter
approximable by $\Pi_T:=\Pi|_{ [0,T]\times{\mathbb R}_+} $,  the restriction of $\Pi$ to the vertical strip $[0,T]\times{\mathbb R}_+$.
We will adopt
a meaningful (and measurable) labeling $(t_i,x_i)$  of  the  atoms  of $\Pi_T$, as  obtained by ranking them
in the  increasing order of $x$-components; this  procedure carries no ambiguity since the ties  occur with zero probability.
The sequence  $(x_i)$ constitutes a univariate Poisson process of intensity $T$, with   concomitant  times $(t_i)$ being an independent sequence of $[0,T]$-uniform marks.
With this enumeration of atoms, the vague convergence to $\Pi_T$ can be more easily described in  terms of  convergence in distribution of the increasing sequence of  uniform order statistics
 \cite[Lemma 3.1]{GJM}.
For instance, the scaled index and  the value of  $(n/T)\,\min (U_1,\ldots,U_n)$ jointly converge  in distribution
to `the least draw'  $(t_1,x_1)$, with independent   $[0,T]$-uniform $t$-component and mean-$T^{-1}$  exponentially distributed $x$-component.

To state  the Poisson process counterpart of Moser's problem we regard
the  generic atom  $(t_i,x_i)$ of $\Pi_T$ as  value $x_i$ drawn at time $t_i$ (but keep in mind that the rank $i$  is not observable at this time).
Formally,  we define a  $\Pi$-adapted {\it stopping time  with strict horizon $T$} to be a random variable $\tau$ whose realisation belongs to  $\{t_i\}$ (thus $\tau<T$ a.s.)
and such that the event $\{\tau\leq t\}$ is measurable relative to $\Pi_t$ for every $t\geq0$.
Thus at time $t_i\in[0,T]$ a decision prescribed by $\tau$  is  made
on base of the information conveyed by the
 draws $(t_j,x_j)$ with $t_j\leq t_i$.
The $x$-value coupled with $\tau$ is the accepted draw value $\xi_\tau$, such that  $\xi_\tau=x_i$ whenever $\tau=t_i$.
For fixed $0<T<\infty$
the optimal stopping problem  amounts to finding the stopping value
\begin{equation}\label{POS}
v(T):=\inf_{\tau<T}\me\,\xi_\tau
\end{equation}
and (if exists) the optimal stopping time $\tau(T)$,
where the infimum is taken over all stopping times with  strict horizon $T$.
We stress that in the problem  (\ref{POS}) there are infinitely many draws observed by any time $t>0$, which differs from
other marked Poisson process models  \cite{BrussSwan, Katriel} like   stopping $\Pi$ restricted to a bounded rectangle.

\paragraph{Notation.} Throughout, an atom $(t_i,x_i)$ is interpreted as the draw value $x_i$ observed at time $t_i$; for a stopping time $\tau$ the accepted draw value is denoted $\xi_\tau$, so that $\xi_\tau=x_i$ on the event $\tau=t_i$.

\paragraph{Threshold stopping times.} Let
 $f:[0,T)\to {\mathbb R}_+$ be a measurable  function, locally bounded on $[0,T)$, with primitive
 $F(t)=\int_0^t f(s)ds$ diverging as  $t\to T$.
Then
\begin{equation}\label{thr}
\tau_f:=\min\{t_i: x_i< f(t_i)\}
\end{equation}
defines a  stopping time with strict horizon $T$. Indeed,
$\Pi_T$ restricted to the subgraph of $f$ is a Poisson point process with mean measure $f(t)dt$ in the time variable, which has infinitely many atoms overall, but finitely many
within each interval $[0,t], 0<t<T$. (Whether the inequality in (\ref{thr}) is strict is immaterial for continuous losses; for the integer-valued losses in Section \ref{S3} the strict form is the convenient one.)
By the same token,
\begin{equation}\label{tau-law}
\prob(\tau_f>t)=e^{-F(t)},~~~ \prob(\tau_f\in dt,\ \xi_{\tau_f}\in dx)=e^{-F(t)}\,1(0\leq x<f(t))\,dt\,dx,
\end{equation}
so  $f$ is the  failure rate of $\tau_f$, and conditionally on $\tau_f=t$ the accepted value is uniform on $[0,f(t)]$.

\paragraph{Problems with finite intensity.} The proofs of optimality below rely on comparison with  a truncated problem in which only finitely many atoms are relevant. The following lemma is a version of the Cayley--Moser problem with Poisson arrival of choice opportunities \cite{Katriel},
it also fits in the framework of Elfving's house-selling problem \cite{CRS}.

\begin{lemma}\label{L1} For $L>0$, let $\ell:[0,L)\to[0,L]$ be a nondecreasing function, and let $T>0$. Consider stopping times $\rho$ adapted to $\Pi$, with values in $\{t_i: t_i<T,\ x_i<L\}\cup\{T\}$, and the loss
$\ell(\xi_\rho)$ if $\rho<T$, and $L$ if $\rho=T$. The minimal expected loss $\varphi(T)$ is the unique solution of
\begin{equation}\label{phiODE}
\varphi'(T)=-\int_0^L\big(\varphi(T)-\ell(x)\big)_+dx,~~~~\varphi(0)=L,
\end{equation}
and it is achieved by the rule which accepts the first atom $(t_i,x_i)$ with $x_i<L$ and $\ell(x_i)\leq \varphi(T-t_i)$, and otherwise stops at $T$.
\end{lemma}

\proof The atoms of $\Pi_T$ below level $L$ arrive at the epochs of a Poisson process of rate $L$, and carry independent $[0,L)$-uniform marks. These atoms are independent of the atoms above level $L$; the latter enter the decision only as an external randomisation, which cannot decrease the minimal expected loss \cite{CRS}. Thus we deal with a problem with almost surely finitely many decision epochs and time-homogeneous dynamics, in which the optimality principle applies in the familiar form. Conditioning on the first arrival gives
$$
\varphi(T)=Le^{-LT}+\int_0^T Le^{-Lr}\left(\frac{1}{L}\int_0^L \min\big(\ell(x),\varphi(T-r)\big)\,dx\right)dr,
$$
whose unique bounded solution is continuously differentiable and satisfies (\ref{phiODE}), as is seen by substituting $r\mapsto T-r$ and differentiating. Since the right-hand side of (\ref{phiODE}) is Lipschitz in $\varphi$, the solution of the initial value problem is unique.
\endpf

The following solution to (\ref{POS}) and connection to Moser's problem exemplify
more general asymptotic results  on optimal stopping for independent  sampling, specifically those concerning source distributions belonging to
 Weibull's  domain of attraction for extrema  \cite{KK1, KK2, Rusch}.
A peculiar feature of (\ref{POS}) is the invariance of the planar Poisson process under the self-similarity transforms $(t,x)\mapsto (\gamma t, x/\gamma), \gamma>0,$ that preserve both intensity and the coordinate-wise
ordering of the atoms.

\begin{theorem}\label{T1} The stopping value and the optimal stopping time for  problem {\rm (\ref{POS})} are given by, respectively,
\begin{eqnarray}\label{re1}
v(T)&=&\frac{2}{T},\\
\label{re2}
\tau(T)&=&\min\left\{t_i:  ~x_i\leq \frac{2}{T-t_i}\right\}.
\end{eqnarray}
\end{theorem}
\noindent
The stopping time $\tau(T)$ is  indeed
well defined  (in particular, has strict horizon $T$) due to  the divergence of the integral of $1/t$ in a vicinity of zero.

\vskip0.3cm

\proof {\it Lower bound.} Fix $L>0$ and let $\tau$ be any stopping time with strict horizon $T$. Consider the truncated problem of Lemma \ref{L1} with $\ell(x)=x$, where only draws below level $L$ may be accepted, but stopping at the terminal time $T$ is allowed and incurs the fixed loss $L$. Define $\tau_L:=\tau$ on the event $\{\xi_\tau<L\}$ and $\tau_L:=T$ otherwise. Then $\tau_L$ is admissible in the truncated problem, and its loss does not exceed $\xi_\tau$. For $\varphi\leq L$ we have $\int_0^L(\varphi-x)_+dx=\varphi^2/2$, so (\ref{phiODE}) becomes $\varphi'=-\varphi^2/2$, $\varphi(0)=L$, whence $\varphi(T)=2L/(TL+2)$. Therefore
$$\me\,\xi_\tau\geq \frac{2L}{TL+2}$$
for every $L$, and letting $L\to\infty$ gives $v(T)\geq 2/T$.

{\it Upper bound.} The stopping time $\tau(T)$ is $\tau_f$ with $f(t)=2/(T-t)$, so  $e^{-F(t)}=(1-t/T)^2$. By (\ref{tau-law}), $\tau(T)$ has density $f(t)e^{-F(t)}$ and the conditional mean of the accepted value given $\tau(T)=t$ is $f(t)/2$, whence
$$
\me\,\xi_{\tau(T)}=\int_0^T \frac{2}{T-t}\left(1-\frac{t}{T}\right)^2\frac{1}{T-t}\,dt=\int_0^T\frac{2}{T^2}\,dt=\frac{2}{T}.
$$
Together with the lower bound this proves (\ref{re1}) and the optimality of (\ref{re2}).
\endpf

\begin{theorem} \label{T2}
For  $\sigma_n$ or any other  sequence of  threshold stopping times in Moser's problem,
\begin{eqnarray*}\label{n-tau-aa}
\tilde{\sigma}_n:=\min\{1\leq j\leq  n: U_j\leq a_{n-j}\},~~~a_0:=\infty,
\end{eqnarray*}
if  the thresholds $(a_n)$ satisfy $n a_n\to 2$, then the convergence in distribution holds,
$$
\left(\frac{T\tilde{\sigma}_n }{n}, \,\frac{nU_{\tilde{\sigma}_n}}{T} \right)\stackrel{d} {\to}
 \left(\tau(T), \xi_{\tau(T)} \right) ~~~({\rm as}~n\to\infty),
$$
together with the convergence of marginal expected values.
\end{theorem}

\proof By the self-similarity, $(\tau(T),\xi_{\tau(T)})\stackrel{d}{=}(T\tau(1),T^{-1}\xi_{\tau(1)})$, hence it is enough to consider $T=1$. Let
$p_j:=\prob(\tilde\sigma_n\geq j)=\prod_{i=1}^{j-1}(1-a_{n-i}\wedge 1)$. By independence of the draws,
$$
\prob(\tilde\sigma_n=j,\ U_{\tilde\sigma_n}\leq y)=p_j\,\min(y,a_{n-j}),~~~0\leq y\leq 1.
$$
Fix $0<s<1$. For $j\leq sn$ all indices $n-i$ in the product exceed $(1-s)n$, so that $a_{n-i}=(2+o(1))/(n-i)$ uniformly, and
$$\log p_j=-\sum_{i=1}^{j-1}\frac{2+o(1)}{n-i}+O\Big(\frac{1}{n}\Big)\to 2\log(1-r),~~~j/n\to r,$$
uniformly in $r\in[0,s]$; also $n a_{n-j}\to 2/(1-r)$. Hence, as a Riemann sum,
$$
\prob\Big(\frac{\tilde\sigma_n}{n}\leq s,\ nU_{\tilde\sigma_n}\leq x\Big)=\sum_{j\leq sn}p_j\min\Big(\frac{x}{n},a_{n-j}\Big)\to \int_0^s(1-r)^2\min\Big(x,\frac{2}{1-r}\Big)dr,
$$
which by (\ref{tau-law}) with $f(r)=2/(1-r)$ is $\prob(\tau(1)\leq s,\ \xi_{\tau(1)}\leq x)$. Since 
$$\prob(\tilde\sigma_n>sn)=p_{\lfloor sn\rfloor+1}\to(1-s)^2,$$
 which is small for $s$ close to $1$, the family of laws is tight and the joint convergence in distribution follows.

Now, $\me\,\tilde\sigma_n/n\to\me\,\tau(1)$ follows by bounded convergence. The convergence of expected values $n\,\me\,U_{\tilde\sigma_n}\to 2=\me\,\xi_{\tau(1)}$ is a special case of \cite[Theorem 3.1]{KK2}. 
To see this, note that \cite{KK2} deals with the {\it maximisation} problem, so we apply their result to $X_j:=-U_j$, with distribution function $F(x)=1+x$ for $-1\leq x\leq 0$. This $F$ belongs to the domain of attraction of $G^\alpha_{\rm III}$ with $\alpha=1$ and $x_F=0$, and the norming constants (\cite{KK2}, Equation (1.6)) are $n$ and $0$. Their threshold stopping time (\cite{KK2}, Equation (1.14)) with thresholds $u_n:=-a_n$ stops at the first $k<n$ with $-U_k>-a_{n-k}$, hence coincides almost surely with $\tilde\sigma_n$. Since $n(1-F(u_n))=n\,(a_n\wedge 1)\to 2$, Theorem 3.1 from that paper applies with $\eta=2>\alpha^{-1}$ and yields
$$\lim_{n\to\infty} n\,\me(-U_{\tilde\sigma_n})=f_{\rm III}(2;1)=-\big(1+1\big)^{-1}\,2^{1+1}\,\big(2-1\big)^{-1}=-2,$$
that is, $n\,\me\, U_{\tilde\sigma_n}\to2$.
\endpf

\subsection{Variational principles}

We sketch three alternative arguments for (\ref{re1}), (\ref{re2}), arranged by their narrowing scope of  applicability,  to explain where the solution comes from and how it generalises.
The first is  just the universal principle of optimality,
the second is apparently not that common (but see  \cite[Remark 5.3]{SK}), and the third is based on the self-similarity previously
employed in the    best-choice and other  sequential selection problems  \cite{GnedinSPA,  GKS, GSe}.

\paragraph{Dynamic programming}
The stopping value $v(T)$ is clearly nonincreasing, and satisfies  $v(0):=v(0+)=\infty$ since the expected value of the lowest  draw of $\Pi_T$ is $T^{-1}$
(prophet's value).
Taking for granted  that  the function is finite  and sufficiently smooth,
an application of  the optimality principle  combined  with the  independence properties of $\Pi_T$ dictate to stop at the earliest draw $(t_i,x_i)$  satisfying $x_i\leq v(T-t_i)$.
The  total probability decomposition over a small time then yields
$$v(T)=\left( \int_0^{v(T)} x\, dx\right)dT+v(T-dT)(1-v(T)dT),$$
which leads to the differential equation
\begin{equation}\label{ODE}
v'(T)=-\frac{1}{2} v^2(T),
\end{equation}
appearing as
an obvious continuous analogue of (\ref{BE}).  The temporal variable $T$ has the meaning of the problem horizon.
Separating the variables in (\ref{ODE}) we readily find that   $v(T)=2/T$ is the      unique finite solution on $(0,\infty)$ which  satisfies  the boundary condition $v(0)=\infty$.

\paragraph{The calculus of variations approach}   Let
 $f:[0,T)\to {\mathbb R}_+$ be a  smooth strictly increasing function with  $F(t)\to\infty$ as $t\to T$, and consider the threshold stopping time (\ref{thr}).
By (\ref{tau-law}),
\begin{equation}\label{mean-f}
\me \xi_{\tau_f}= \frac{1}{2}\int_0^T e^{-F(t)}f^2(t)dt.
\end{equation}
For $0\leq x\leq f(0)$
\begin{equation}\label{dis-tau1}
\prob(\xi_{\tau_f}\leq x)= x\,\me \tau_f=x\,\int_0^T e^{-F(t)}dt,
\end{equation}
hence on this interval the distribution  of the stopped draw value is uniform.
For $x>f(0)$ we can express the survival function as
\begin{equation}\label{dis-tau2}
\prob(\xi_{\tau_f}> x)=\int_{f^{-1}(x)}^T e^{-F(t)}(f(t)-x)dt.
\end{equation}

To minimise  the functional  (\ref{mean-f}) in $F$ write
the Euler-Lagrange equation
$$
\frac{d}{dt}\frac{d}{df}\left( \frac{1}{2}e^{-F}f^2\right)=\frac{d}{dF}\left( \frac{1}{2}e^{-F}f^2\right).
$$
 With the account of $F'=f$ this reduces to the first order differential equation
\begin{equation}\label{odef}
f'(t)=\frac{1}{2}f^2(t),
\end{equation}
whose  unique  solution  finite on $[0,T)$ with $f(T)=\infty$ is
\begin{equation}\label{beta2}
f(t)=\frac{2}{T-t}.
\end{equation}
Passing to the time reversal of $f$ we are back to  the Wald-Bellman equation (\ref{ODE}) for the stopping value  $v(t)=f(T-t)$.

\paragraph{Self-similarity}
 Generalising (\ref{beta2}),
we call  {\it beta strategy}, denoted $\beta_b(T)$, a stopping time (\ref{thr}) with the hyperbolic threshold  function $f(t)=b/(T-t)$.
For every $b>0$ this is  a well defined stopping time with strict horizon $T$, however $\me \,\xi_{\beta_{b}(T)}<\infty$ requires $b>1$.

The self-similarity of $\Pi$ makes
beta strategies inherent to the planar Poisson framework.
This entails the homogeneity of the stopping value, $v(aT)=a^{-1}v(T)$,  which gives $v(T)=v(1)/T$ and suggests that the beta strategy with $b=v(1)$ is optimal for every $T$.
The  problem
(\ref{POS}) is thus  reduced to  finding $v(1)$ by means of optimisation over  the  one-parameter class of beta strategies.

Using $\beta_{b}(T)\stackrel{d}{=}T\,\beta_b(1)$ we resort
 to $T=1$, to  first find that $\beta_{b}(1)$ has a beta distribution:
\begin{equation}\label{beta-d}
{\mathbb P}(\beta_{b}(1)>t)=\exp\left(-\int_0^t \frac{b\,ds}{1-s}\right)=(1-t)^b, ~~t\in[0,1].
\end{equation}
Now (\ref{mean-f}) specialises as
\begin{equation}\label{mean-b}
\me\,\xi_{\beta_{b}(1)}=\frac{1}{2}\int_0^1 (1-t)^b \left(\frac{b}{(1-t)}\right)^2  dt=\frac{b^2}{2(b-1)},
\end{equation}
which is minimised  at $b=2$.
The optimality of  $\beta_{2}(T)$  follows for every $T>0$.

The minimiser $b=2$  could have been  determined straight from the equation
$$
b=\frac{b^2}{2(b-1)},
$$
saying that if a  draw occurs at time $0$ and hits the threshold, the conditional continuation and stopping risks are  equal.
This equation exemplifies  a  well known condition called  differently as `continuous fit principle'  \cite{Peskir}  or sometimes `balance at the boundary' \cite{SK}.

For the original  discrete-time  Moser's problem the selfsimilarity does not hold  exactly.
But for large $n$  the property
shows up in the extreme-value range of draws in the form of
   asymptotic homogeneity
$ V_{\lfloor T n \rfloor}\sim T^{-1}V_n$.

\subsection{Details on beta strategies}

By the selfsimilarity the stopped  draw satisfies
$$
(\beta_b(T),\xi_{\beta_{b}(T)})\stackrel{d}{=} (T\beta_b(1),T^{-1}\xi_{\beta_{b}(1)}),
$$
thus  we lose no generality by setting $T=1$. We write $\beta_b$  as a shorthand for  $\beta_b(1)$, and assume  that $b>1$.
According to (\ref{beta-d}), $\beta_b$ follows beta$(1,b)$ distribution, which has moments
$$\me \beta_b=\frac{1}{b+1},~{\rm Var}(\beta_b)= \frac{b}{(b+1)^2(b+2)}.$$

In the case $b=2$ we are back to  (\ref{MM}). The asymptotics follow from the convergence in distribution (cf Theorem \ref{T2})
taken together with the uniform integrability of the scaled stopping time  $\sigma_n/n\leq 1$.

The distribution of the stopped draw value is much more intriguing.
Explicit integration in (\ref{dis-tau1}), (\ref{dis-tau2}) outputs  a piecewise density
\begin{equation}\label{un-par}
\prob(\xi_{\beta_b}\in dx)=\begin{cases} \dfrac{1}{b+1}\,, ~~~~~~~~~~~~0\leq x\leq b,\\
\dfrac{1}{b+1}\left(\dfrac{b}{x}\right)^{b+1}\,,~~~x>b.
\end{cases}
\end{equation}
The distribution  appeared in
\cite{KK1} (the last case in Equation (2.11), with $\alpha=1$), see also \cite[Equation 4.13]{SK}.
We may regard  distribution (\ref{un-par}) as  a mixture of  the $[0,b]$-uniform  and the Pareto distributions, the latter  with both shape and scale parameters equal to $b$.
In particular, the moments are finite only for orders less than $b$.  
Unlike the minimal order statistic approximable by the exponential distribution,
the stopped draw in Moser's problem $(b=2)$  is asymptotically  heavy-tailed: the density equals $8/(3x^3)$ for $x>2$, so the mean is finite but the second moment is infinite.

\section{Sampling from  the discrete uniform distribution}\label{S3}

\subsection{The discrete time stopping problem  } \label{3.1}

We prefer to consider the discrete uniform distribution on $\{0,1,\ldots,N-1\}$, that is with  $0$ being the least possible draw value,
to match neatly with the support of the continuous  uniform distribution.
The discrete-range independent sampling
is then represented
via rounded {\it down} continuous uniform random variables $\lfloor  N U_1\rfloor, \lfloor N U_2\rfloor,\ldots$.
The optimal stopping problem with $n$ draws requires  to find the stopping value
\begin{equation}\label{disMP}
V_{n,N}:=\inf_{\sigma\leq n} {\mathbb E}\,\lfloor  N   U_\sigma \rfloor,
\end{equation}
where  the infimum is taken over the same class of stopping times adapted to $U_j$'s  as in (\ref{ContiMP}).
Having  $U_j$'s (in excess of  $\lfloor NU_j\rfloor$'s)  observable  brings no advantage to the online decision maker,
as  the  information contained in independent fractional parts $NU_j-\lfloor NU_j\rfloor$ is just a superfluous randomisation
that cannot help improving the stopping value.

The optimality equation becomes
\begin{equation}\label{discV}
V_{n,N}-V_{n-1,N}=-\frac{1}{N}\sum_{k=0}^{N-1}(V_{n-1,N}-k)_+, ~~~n\geq 1,~~~ V_{0,N}:=N-1.
\end{equation}
Comparing with (\ref{BE}),
this is more difficult to handle since the critical thresholds now depend on both $n$ and $N$.
The  numerical values computed from (\ref{discV}) and shown with 2 decimal places demonstrate the opposite  directions of  monotonicity { (rows $N$, columns $n$)}:
$$
\left(V_{n,N}\right)_{2\leq N\leq10, ~1\leq n\leq 10}=
\begin{bmatrix}
0.50 & 0.25 & 0.13 & 0.06 & 0.03 & 0.02 & 0.01 & 0.00 & 0.00 & 0.00 \\
1.00 & 0.67 & 0.44 & 0.30 & 0.20 & 0.13 & 0.09 & 0.06 & 0.04 & 0.03 \\
1.50 & 1.00 & 0.75 & 0.56 & 0.42 & 0.32 & 0.24 & 0.18 & 0.13 & 0.10 \\
2.00 & 1.40 & 1.04 & 0.82 & 0.66 & 0.53 & 0.42 & 0.34 & 0.27 & 0.22 \\
2.50 & 1.75 & 1.33 & 1.06 & 0.87 & 0.73 & 0.60 & 0.50 & 0.42 & 0.35 \\
3.00 & 2.14 & 1.65 & 1.32 & 1.09 & 0.92 & 0.79 & 0.68 & 0.58 & 0.50 \\
3.50 & 2.50 & 1.94 & 1.58 & 1.31 & 1.11 & 0.95 & 0.84 & 0.73 & 0.64 \\
4.00 & 2.89 & 2.26 & 1.84 & 1.54 & 1.31 & 1.13 & 0.99 & 0.88 & 0.78 \\
4.50 & 3.25 & 2.55 & 2.09 & 1.76 & 1.51 & 1.31 & 1.14 & 1.02 & 0.91
\end{bmatrix}
$$
Easily enough, as in Moser's problem, $V_{n,N}$ decreases with the number of trials $n$.
The increasing pattern in $N$, though intuitive  since  the distribution puts less weight on smaller values, is not immediate but can be inferred with the aid of the
coupling
with continuous uniform draws. Indeed, we have  $\me \lfloor NU_\sigma\rfloor$  nondecreasing in $N$ for any
 $(U_j)$-adapted  stopping time $\sigma\leq n$, therefore the strategy optimal for  the $\lfloor NU_j\rfloor$'s  will act only    suboptimally
when applied to $\lfloor (N+1)U_j\rfloor$'s.

We are interested in the asymptotic regime  $n,N\to\infty, n/N\to T$   with $0<T<\infty$ kept as parameter.
By the virtue of  $\lfloor NU_j\rfloor< NU_j<\lfloor NU_j\rfloor+1$ (holding almost surely),
coupling  with  Moser's problem (\ref{ContiMP})  implies
\begin{equation}\label{VNcoupling}
                              V_{n,N}< NV_n<V_{n,N}+1,
\end{equation}
where $V_{n,N}+1$ is the stopping value for sampling from the  (more common)  uniform distribution on $\{1,2,\ldots,N\}$.
These elementary bounds allow us to squeeze the asymptotic stopping value as
\begin{equation}\label{CompU}
\left( \frac{2}{T}-1\right)_+\leq \liminf_{n/N\to T} V_{n,N}\leq \limsup_{n/N\to T} V_{n,N}\leq  \frac{2}{T},
\end{equation}
and suggest to seek for an analogue of (\ref{Mos}). The lower bound will be improved to $2/T-1/2$ in Proposition \ref{P1}.

\subsection{The Poisson counterpart}\label{S3.2}

For $n/N\to T$ the point scatter $\{(j/N, NU_j), 1\leq j\leq n\}$ has the same point-process limit $\Pi_T$ as under scaling (\ref{SP}).
The number of $\Pi_T$-atoms with $x$-values falling in a {\it box} $[k,k+1)$ before time $T$ has the Poisson distribution with mean $T$, and the counts for boxes $k=0,1,\ldots$ are  independent.
In particular,  each box remains empty at time $T$ with probability $e^{-T}$.
This implies  that  the least rounded draw value $\lfloor x_1\rfloor$  has  the geometric distribution, whence  the asymptotic prophet's benchmark
$$\me \min(\lfloor NU_1\rfloor,\ldots, \lfloor NU_n\rfloor)\to\frac{1}{e^T-1},$$
which improves upon the lower bound from  (\ref{CompU})
for $T>1.594$, and upon the improved lower bound $2/T-1/2$ of Proposition \ref{P1} for $T>3.830$.

In the limit dynamical picture,  the draws from different boxes occur at epochs of  independent unit-rate Poisson processes on $[0, T]$. Stopping with a draw  from box $[k,k+1)$
incurs loss $k$. To frame  this formally,  we pose  the Poisson process counterpart of  (\ref{disMP}) as the minimisation problem of finding the stopping value
\begin{equation}\label{dePo}
u(T):=\inf_{\tau<T} \me\lfloor\xi_\tau\rfloor,
\end{equation}
where  the infimum  expands over the same class of   stopping times as in (\ref{POS}).

We pause to give a definition for the rest of the paper. Under {\it loss structure} in a stopping problem related to a planar Poisson process we shall understand a function in the variables $(t,x)$,
to evaluate the outcome of stopping at this location should there be an atom of the process.  Thus for (\ref{dePo}) the loss structure is $(t,x)\mapsto \lfloor x\rfloor$.

Since $\lfloor x\rfloor\leq x<\lfloor x\rfloor+1$, comparison with Theorem \ref{T1} shows
\begin{equation}\label{uBounds0}
\left(\frac{2}{T}-1\right)_+\leq u(T)\leq \frac{2}{T}.
\end{equation}
We introduce a function
$$\Phi(z):=\sum_{j=0}^\infty (z-j)_+,~~~z\geq 0,$$
which is
 convex, strictly increasing and piecewise linear with slope $j+1$ on the interval $(j,j+1)$.

\paragraph{The optimality equation.}
The stopping value function (\ref{dePo}) is expected to satisfy the optimality equation
\begin{equation}\label{ucont}
u'(T)=-\Phi(u(T))=-\sum_{j=0}^\infty \big(u(T)-j\big)_+, ~~~u(0)=\infty,
\end{equation}
which is analogous to (\ref{discV}).
We will first solve (\ref{ucont}) and then prove in Theorem \ref{T3} that the unique solution is indeed the stopping value.
A continuous solution, strictly decreasing to $0$ as $T\to\infty$, defines, for every $k\geq 1$, a unique {\it cutoff} $\delta_k$ satisfying
\begin{equation}\label{fit}
u(\delta_k)=k,
\end{equation}
and let $\delta_0:=\infty$. The sequence of cutoffs is  strictly decreasing to $0$ as $k\to\infty$,
which altogether implies for  $k\geq0$ that
$$
k<u(T)<k+1\Longleftrightarrow \delta_{k+1}<T<\delta_k.
$$

The cutoff  form of strategy is quite common in sequential decision problems with discrete observables (see e.g. \cite{Gianini,  GS1,GD1, GD2,  Mucci1, Mucci2}).

In the range  $\delta_{k+1}<T<\delta_{k}$ the sum in (\ref{ucont}) has positive terms with  $0\leq j\leq k$ only, and the equation simplifies
\begin{equation}\label{linODE}
u'(T)=-(k+1)u(T)+  \frac{k(k+1)}{2},
\end{equation}
with the general solution being
\begin{equation}\label{genU}
u(T)=\frac{k}{2}+e^{-(k+1)T}C.
\end{equation}
Invoking  condition (\ref{fit})  at the right endpoint gives a piecewise solution
\begin{equation}\label{par-sol}
u(T)=\begin{dcases} \tfrac{k}{2}\left( e^{(k+1)(\delta_{k}-T)}+1\right),~~\delta_{k+1}\leq T\leq \delta_{k}, ~k\geq 1,\\
e^{\delta_1-T}, ~~~~~~~~~~~~~~~~~~~~~~T\geq \delta_1.
\end{dcases}
\end{equation}
Matching  with (\ref{fit})  at  the left endpoint
 results in a recursion for the cutoffs 
\begin{equation}\label{reccut}
e^{\delta_{k}-\delta_{k+1}}=\left(1+\frac{2}{k}\right)^{\tfrac{1}{k+1}},~~~k\geq1,
\end{equation}
that implies
$$
e^{\delta_1-\delta_{k}}=\prod_{j=1}^{k-1}\left(1+\frac{2}{j}\right)^{\tfrac{1}{j+1}}.
$$
Letting $k\to\infty$ we should have $e^{\delta_k}\to 1$, thus first for $k=1$ then  for all $k\geq 1$ we obtain the infinite product formula
\begin{equation}\label{edelta}
e^{\delta_k}=\prod_{j=k}^\infty  \left(1+\frac{2}{j}\right)^{\tfrac{1}{j+1}}, ~~~k\geq1.
\end{equation}
Conversely, defining the cutoffs by (\ref{edelta}) and $u$ by (\ref{par-sol}) one obtains a continuous, strictly decreasing function on $(0,\infty)$ with $u(0+)=\infty$, which satisfies (\ref{ucont}) on each interval $(\delta_{k+1},\delta_k)$ and is continuously differentiable also at the cutoffs, since both one-sided derivatives at $\delta_k$ are equal to $-\Phi(k)=-k(k+1)/2$.
This solution is unique. Indeed, since (\ref{ucont}) is autonomous and $\Phi$ is locally Lipschitz, every solution with values in $(0,\infty)$ is a time shift $T\mapsto u(T+c)$ of the solution just constructed (choose $c$ to match the value at any one point), and the boundary condition $u(0)=\infty$ forces $c=0$. In the same way, $v(T)=2/T$ is the only solution of (\ref{ODE}) exploding at $0$. Thus the seemingly improper initial condition $u(0)=\infty$ singles out the right solution; in the proof of Theorem \ref{T3} it appears as the limit of the proper initial conditions $\varphi(0)=L$ of Lemma \ref{L1}.

\begin{lemma}\label{L2} Let $u$ be any solution of $u'=-\Phi(u)$ on $(0,\infty)$ with $u(0+)=\infty$. Then, for all $T>0$,
\begin{equation}\label{uBounds}
\frac{2}{T}-\frac{1}{2}\ \leq\ u(T)\ \leq\ \min\left(\frac{2}{T},\ \frac{2}{T(1-T^2/48)}-\frac12\right),
\end{equation}
where the second expression in the minimum is used for $T<\sqrt{48}$.
\end{lemma}

\proof Writing $z=k+r$ with $k=\lfloor z\rfloor$ and $0\leq r<1$, we have $\Phi(z)=(k+1)z-k(k+1)/2$, which rearranges as
$$
\Phi(z)=\frac{1}{2}\Big(z+\frac12\Big)^2-\theta(z),~~~~\theta(z):=\frac{1}{8}-\frac{r(1-r)}{2}\in\Big[0,\frac{1}{8}\Big].
$$
In particular, $\Phi(z)\geq z^2/2$, whence $(1/u)'=\Phi(u)/u^2\geq 1/2$ and $u(T)\leq 2/T$. Similarly, $Z:=u+\frac12$ satisfies $Z'=-Z^2/2+\theta(Z-1/2)$, so that
$$\Big(\frac{1}{Z}\Big)'=\frac{1}{2}-\frac{\theta(Z-1/2)}{Z^2}.$$
Since $1/Z(0+)=0$, integration gives $1/Z(T)\leq T/2$, whence the lower bound in (\ref{uBounds}). Now using $Z(t)\geq 2/t$ we integrate to obtain

$$\frac{1}{Z(T)}\geq \frac{T}{2}-\frac{1}{8}\int_0^T\frac{t^2}{4}dt=\frac{T}{2}\left(1-\frac{T^2}{48}\right),$$ 
which gives the remaining part of the upper bound.
\endpf

The bounds imply  $u(T)=2/T-1/2+O(T)$ as $T\to 0$. Applied at the cutoffs, where $u(\delta_k)=k$, they yield the tight two-sided estimate
\begin{equation}\label{cutBounds}
\frac{2}{k+1/2}\leq \delta_k,~~~~~~\delta_k\Big(1-\frac{\delta_k^2}{48}\Big)\leq \frac{2}{k+1/2},
\end{equation}
so  $\delta_k=2/(k+1/2)+O(k^{-3})$ for large $k$. We will also use the generalised inverse of $k\mapsto\delta_k$,
$$\varkappa(s):=\max\{k\geq 0:\ \delta_k\geq s\},~~~s>0,$$
which satisfies $\varkappa(s)=\lfloor u(s)\rfloor$ off the cutoffs; by (\ref{uBounds}),
\begin{equation}\label{kappaBounds}
\frac{2}{s}-\frac12\leq \varkappa(s)+1\leq \frac{2}{s}+1.
\end{equation}

\begin{theorem}\label{T3}
The stopping value function {\rm(\ref{dePo})} is given by {\rm (\ref{par-sol})} with the cutoffs
\begin{equation}\label{delta-cut}
\delta_k=\sum_{j=k}^\infty  \frac{\log\left(1+{2}/{j}\right)}{j+1}.
\end{equation}
The optimal stopping time for the problem {\rm (\ref{dePo})} with loss structure $(t,x)\mapsto\lfloor x\rfloor $  is
\begin{equation}\label{nu}
\nu(T)=\min\{t_i:   t_i\geq (T-\delta_k)_+ {\rm ~for~}k=\lfloor x_i\rfloor\}.
\end{equation}
\end{theorem}

\proof Let $u$ be the solution (\ref{par-sol}) of (\ref{ucont}) with cutoffs (\ref{delta-cut}), unique by the above; it is continuous, strictly decreasing from $u(0+)=\infty$ to $0$, and satisfies the bounds of Lemma \ref{L2}. For the moment, denote by $\bar u(T)$ the stopping value (\ref{dePo}); we need to show that $\bar u=u$ and that $\nu(T)$ achieves $u(T)$.

{\it The truncated problem.} Let $L\geq1$ be an integer, and let $c_L>0$ be defined by $u(c_L)=L$; then $c_L\to0$ as $L\to\infty$. Consider the setting of Lemma \ref{L1} with $\ell(x)=\lfloor x\rfloor$ and terminal loss $L$. For $0\leq\varphi\leq L$ we have $\int_0^L(\varphi-\lfloor x\rfloor)_+dx=\Phi(\varphi)$, so (\ref{phiODE}) becomes $\varphi'=-\Phi(\varphi)$, $\varphi(0)=L$. By the autonomy of the equation, its solution is the time shift
$$\varphi_L(T)=u(T+c_L).$$

{\it Lower bound.} Let $\tau$ be a stopping time with strict horizon $T$, and let $\tau_L:=\tau$ if $\xi_\tau<L$, and $\tau_L:=T$ otherwise. Then $\tau_L$ is admissible in the truncated problem, and its loss does not exceed $\lfloor\xi_\tau\rfloor$ (on the event $\xi_\tau\geq L$ we have $\lfloor\xi_\tau\rfloor\geq L$). Hence $\me\lfloor\xi_\tau\rfloor\geq\varphi_L(T)=u(T+c_L)$, and letting $L\to\infty$ gives $\bar u(T)\geq u(T)$.

{\it Upper bound.} Since $u$ is decreasing with $u(\delta_k)=k$, the condition $t_i\geq(T-\delta_k)_+$ in (\ref{nu}) is equivalent to $\lfloor x_i\rfloor\leq u(T-t_i)$, that is to $x_i<\varkappa(T-t_i)+1$. Thus $\nu(T)=\tau_f$ is a threshold stopping time (\ref{thr}) with the step function $f(t)=\varkappa(T-t)+1$. By (\ref{kappaBounds}), $f(t)\leq 2/(T-t)+1$ and
$$F(t)=\int_0^t f(s)\,ds\,\geq\,\int_0^t\Big(\frac{2}{T-s}-\frac12\Big)ds=2\log\frac{T}{T-t}-\frac{t}{2}\,\geq\, 2\log\frac{T}{T-t}-\frac{T}{2},~~~0\leq t<T,$$
whence
\begin{equation}\label{nu-est}
f(t)\leq\frac{2}{T-t}+1,~~~~e^{-F(t)}\leq e^{T/2}\Big(\frac{T-t}{T}\Big)^2;
\end{equation}
in particular $F(t)\to\infty$ as $t\to T$, so $\nu(T)$ has strict horizon $T$.

Now fix $L$ so large that $c:=c_L<T$, and compare $\nu(T)$ with the optimal rule $\rho_L$ for the truncated problem with horizon $T-c$. By Lemma \ref{L1}, $\rho_L$ accepts the first atom with $t_i<T-c$, $x_i<L$ and $\lfloor x_i\rfloor\leq\varphi_L(T-c-t_i)=u(T-t_i)$, and otherwise stops at $T-c$ with loss $L$. For $t_i<T-c$ we have $u(T-t_i)<u(c)=L$, so the condition $x_i<L$ is implied by the last one. Therefore $\rho_L=\nu(T)$ on the event $\{\nu(T)<T-c\}$, and $\rho_L=T-c$ otherwise, whence
$$
\me\lfloor\xi_{\nu(T)}\rfloor=\varphi_L(T-c)+\me\Big[\big(\lfloor\xi_{\nu(T)}\rfloor-L\big)1(\nu(T)\geq T-c)\Big]\leq u(T)+\me\big[\xi_{\nu(T)}1(\nu(T)\geq T-c)\big].
$$
By (\ref{tau-law}) and (\ref{nu-est}), the last expectation is
$$
\int_{T-c}^T e^{-F(t)}\frac{f^2(t)}{2}\,dt\leq \frac{e^{T/2}}{2T^2}\int_{T-c}^T\big(2+(T-t)\big)^2dt\leq \frac{e^{T/2}(2+T)^2}{2T^2}\,c,
$$
which tends to $0$ as $L\to\infty$. Therefore $\me\lfloor\xi_{\nu(T)}\rfloor\leq u(T)$.

Combining the bounds, $\bar u(T)\leq\me\lfloor\xi_{\nu(T)}\rfloor\leq u(T)\leq \bar u(T)$, which proves both claims of the theorem.
\endpf

Thus the optimal stopping rule $\nu(T)$ is willing to accept any draw value $k\leq x_i<k+1$
(i.e. $\lfloor x_i\rfloor=k$)
starting from the time $(T-\delta_k)_+$, for $k\geq0$, and Lemma \ref{L2} applies to the stopping value $u$.
\paragraph{Smooth fit in the time domain.}
The function $u(T)$ is continuously differentiable, also at the cutoffs $\delta_k$. This is a property of the optimal cutoffs only.
Indeed, consider a strategy which accepts draw value $k$ when the remaining horizon is less than some $\delta_k'$, and let $W$ be its continuation value as a function of the horizon. For horizons slightly above and slightly below $\delta'_k$, $W$ satisfies the linear equations
$$W'=-kW+\frac{(k-1)k}{2}~~~{\rm and}~~~W'=-(k+1)W+\frac{k(k+1)}{2},$$
respectively, so the one-sided derivatives of $W$ at $\delta_k'$ differ by $W(\delta'_k)-k$. Therefore $W$ is differentiable at $\delta'_k$ if and only if $W(\delta'_k)=k$, that is, if at the cutoff the loss incurred by accepting draw value $k$ equals the continuation value. In other words, smooth fit in the time variable is equivalent to the continuous fit (`balance at the boundary') in the draw-value variable, which is the condition (\ref{fit}) we used to determine the cutoffs.

This should be distinguished from the usual smooth fit principle in the state space \cite{Peskir}. In our setting there is no state variable besides the remaining horizon: the draws arrive as atoms of the Poisson process independently of the past, and the loss $\lfloor x\rfloor$ is discontinuous in $x$, so the function $(T,x)\mapsto\min(\lfloor x\rfloor, u(T))$ has no fit of any kind in $x$. The only free boundary is the set of cutoffs in the time domain, which the deterministic time variable crosses continuously.
In Poisson-driven problems with a genuine state variable the picture is more complex; for instance, in the adaptive Poisson disorder problem the value function satisfies the smooth fit on a part of the optimal stopping boundary and fails it on the rest \cite{BDK}.
A time-domain balance condition of the same nature as (\ref{fit}) determines the optimal boundary in the memoryless best-choice problem for the planar Poisson process \cite{GC}.

A few cutoff values (three decimals) are shown in the table
$$
\begin{array}{ ccc ccc ccc  ccc}
k & 1&2&3   & 4&5&6  &7 &8&9 &10&11   \\
\delta_k & 1.353 & 0.804 &0.573&0.445   &0.364 &0.308 & 0.267 &0.235 & 0.211 &0.191& 0.174\\
\end{array}
$$
The bounds (\ref{cutBounds})  give an excellent approximation: $2/(k+1/2)$ equals $1.333, 0.800, 0.571, 0.444$ for $k=1,2,3,4$, and differs from $\delta_k$ by less than $0.0004$ for $k\geq 5$.
Thus for $k$ large the halved cutoffs are well approximable by $1/(k+1/2)$.

\subsection{Convergence of the discrete problem}

\begin{theorem}\label{T4} As $n,N\to\infty$ and $n/N\to T\in(0,\infty)$,
$$V_{n,N}\to u(T).$$
Moreover, the convergence is uniform for $n/N$ in any compact subset of $(0,\infty)$.
\end{theorem}

\proof The key observation is that the optimality recursion (\ref{discV}) coincides with the Euler scheme with step $h:=1/N$ for the differential equation (\ref{ucont}). Indeed, $V_{m,N}\leq V_{0,N}<N$, so the terms with $k\geq N$ in $\Phi(V_{m-1,N})$ vanish and (\ref{discV}) reads as
$$y_m=y_{m-1}-h\,\Phi(y_{m-1}),~~~~ y_m:=V_{m,N},~~m\geq 1.$$
Let $z_m:=u(mh)$. Fix $0<\varepsilon<T$. On $[\varepsilon,\infty)$ the derivative $u'=-\Phi(u)$ is Lipschitz with constant $K_\varepsilon:=(u(\varepsilon)+1)\Phi(u(\varepsilon))$, because $\Phi$ has slope at most $u(\varepsilon)+1$ on $[0,u(\varepsilon)]$. Hence for $mh\geq\varepsilon+h$
$$z_m=z_{m-1}-h\,\Phi(z_{m-1})+r_m,~~~|r_m|\leq K_\varepsilon h^2.$$
The approximation errors $e_m:=y_m-z_m$ satisfy $e_m=(1-h\lambda_m)e_{m-1}-r_m$, where
$$\lambda_m:=\frac{\Phi(y_{m-1})-\Phi(z_{m-1})}{y_{m-1}-z_{m-1}}$$
(and $\lambda_m:=0$ if $y_{m-1}=z_{m-1}$). Recall that $\Phi$ is convex and piecewise linear with slope $\lfloor z\rfloor+1$ on each interval $(\lfloor z\rfloor,\lfloor z\rfloor+1)$. A difference quotient of such a function lies between its slopes at the endpoints, therefore with $\underline{z}:=\min(y_{m-1},z_{m-1})$ and $\overline{z}:=\max(y_{m-1},z_{m-1})$
$$\underline{z}\,\leq\,\lfloor \underline{z}\rfloor+1\,\leq\, \lambda_m\,\leq\, \lfloor \overline{z}\rfloor+1\,\leq\, \overline{z}+1 .$$

From (\ref{MoserBounds}), (\ref{VNcoupling}) and (\ref{uBounds}),
\begin{equation}\label{yzBounds}
\frac{2}{h(m+H_m+2)}-1< y_m<\frac{2}{hm},~~~~~\frac{2}{hm}-\frac12\leq z_m\leq \frac{2}{hm}.
\end{equation}
Let $m_0:=\lceil \varepsilon N\rceil$. For $m>m_0$ and $N$ large enough, (\ref{yzBounds}) gives $h\lambda_m\leq h(2/\varepsilon+1)<1$, hence
$$
0\leq 1-h\lambda_m\leq 1-h\min(y_{m-1},z_{m-1})\leq \exp\Big(-\frac{2}{m+H_m+1}+h\Big).
$$
Iterating the error recursion from $m_0$ to $n$,
$$
|e_n|\leq |e_{m_0}|\prod_{m=m_0+1}^n\exp\Big(-\frac{2}{m+H_m+1}+h\Big)+\sum_{m=m_0+1}^n|r_m|.
$$
As $N\to\infty$ with $n/N\to T$: the product converges to $(\varepsilon/T)^2e^{T-\varepsilon}$ because $H_m\sim \log m$; the sum is at most $nK_\varepsilon h^2\to0$; and by (\ref{yzBounds}), $\limsup|e_{m_0}|\leq 1$. Consequently
$$\limsup_{n/N\to T}|V_{n,N}-u(n/N)|\leq \Big(\frac{\varepsilon}{T}\Big)^2e^{T},$$
and since $\varepsilon$ is arbitrary and $u$ is continuous, $V_{n,N}\to u(T)$. All estimates are uniform for $T$ in a compact subset of $(0,\infty)$.
\endpf

\begin{rem}\rm Since $u$ is strictly decreasing and the convergence in Theorem \ref{T4} is locally uniform, the cutoffs of the discrete problem also converge. By (\ref{discV}) the optimal rule accepts the draw value $k$ when $m$ further draws remain if and only if $k\leq V_{m,N}$; hence the largest such $m$, divided by $N$, converges to $\delta_k$ for every $k\geq1$. Numerically, $V_{N,N}=1.4958,\ 1.5110,\ 1.5129$ for $N=10^3,10^4,10^5$, compared with $u(1)=1.5132$.
\end{rem}

\paragraph{Example} In the $n=N$ case, that is
for  $n$ draws with replacement from $\{0,1,\ldots,n-1\}$,
the optimal stopping rule dictates  to accept  draw values $0$ and $1$ throughout, $2$ when at most about  $0.804 \,n$ draws remain, $3$ when at most  about $0.573 \,n$ draws remain, etc.
The expected value $V_{n,n}$ with this strategy approaches  $u(1)=1.513\ldots$ (Theorem \ref{T4}).
(The limit   $2.513\ldots$ in Introduction is valid
for sampling  from  the uniform distribution on $\{1,2,\ldots,n\}$.)

The  diagonal  entries  in the table of Section \ref{3.1} and further  computed  values
$$(V_{n,n})_{11\leq n\leq 17}=( 0.95, 0.97, 0.99,1.02,  1.04, 1.05, 1.07)$$
suggest that the sequence  $(V_{n,n})$ is monotonically increasing.
\vskip0.2cm

\begin{figure}[htbp]
\centering
\begin{tikzpicture}
\begin{axis}[width=0.8\textwidth,height=0.5\textwidth,xmin=0.15,xmax=1.6,ymin=0,ymax=12,
 axis lines=left, xlabel={$T$}, xlabel style={at={(axis description cs:1.02,0)},anchor=west},
 ytick={0,2,...,12}, xtick={0.2,0.4,...,1.6}, tick label style={font=\small},
 every axis plot/.append style={line width=0.9pt}, clip=false]
\addplot[gray!70!black, line width=0.7pt, const plot, mark=none] coordinates {(0.1600,12) (0.1739,11) (0.1905,10) (0.2106,9) (0.2354,8) (0.2668,7) (0.3079,6) (0.3640,5) (0.4451,4) (0.5728,3) (0.8038,2) (1.3531,1) (1.6000,1)};
\addplot[black, dashed] coordinates {(0.1682,11.8927) (0.1718,11.6411) (0.1754,11.4000) (0.1791,11.1686) (0.1827,10.9465) (0.1863,10.7330) (0.1900,10.5277) (0.1936,10.3301) (0.1972,10.1398) (0.2009,9.9563) (0.2045,9.7794) (0.2081,9.6087) (0.2118,9.4438) (0.2154,9.2845) (0.2190,9.1304) (0.2227,8.9814) (0.2263,8.8372) (0.2299,8.6975) (0.2336,8.5622) (0.2372,8.4311) (0.2409,8.3039) (0.2445,8.1804) (0.2481,8.0606) (0.2518,7.9443) (0.2554,7.8312) (0.2590,7.7213) (0.2627,7.6145) (0.2663,7.5106) (0.2699,7.4095) (0.2736,7.3110) (0.2772,7.2152) (0.2808,7.1218) (0.2845,7.0308) (0.2881,6.9421) (0.2917,6.8557) (0.2954,6.7713) (0.2990,6.6890) (0.3026,6.6087) (0.3063,6.5303) (0.3099,6.4537) (0.3135,6.3789) (0.3172,6.3058) (0.3208,6.2344) (0.3244,6.1645) (0.3281,6.0963) (0.3317,6.0295) (0.3353,5.9641) (0.3390,5.9002) (0.3426,5.8376) (0.3462,5.7763) (0.3499,5.7163) (0.3535,5.6576) (0.3571,5.6000) (0.3608,5.5436) (0.3644,5.4883) (0.3680,5.4341) (0.3717,5.3810) (0.3753,5.3289) (0.3789,5.2778) (0.3826,5.2276) (0.3862,5.1785) (0.3898,5.1302) (0.3935,5.0828) (0.3971,5.0363) (0.4008,4.9906) (0.4044,4.9458) (0.4080,4.9017) (0.4117,4.8584) (0.4153,4.8159) (0.4189,4.7742) (0.4226,4.7331) (0.4262,4.6927) (0.4298,4.6531) (0.4335,4.6141) (0.4371,4.5757) (0.4407,4.5380) (0.4444,4.5008) (0.4480,4.4643) (0.4516,4.4284) (0.4553,4.3931) (0.4589,4.3583) (0.4625,4.3240) (0.4662,4.2903) (0.4698,4.2571) (0.4734,4.2245) (0.4771,4.1923) (0.4807,4.1606) (0.4843,4.1294) (0.4880,4.0986) (0.4916,4.0683) (0.4952,4.0385) (0.4989,4.0090) (0.5025,3.9800) (0.5061,3.9515) (0.5098,3.9233) (0.5134,3.8955) (0.5170,3.8682) (0.5207,3.8412) (0.5243,3.8145) (0.5279,3.7883) (0.5316,3.7624) (0.5352,3.7368) (0.5388,3.7116) (0.5425,3.6868) (0.5461,3.6622) (0.5497,3.6380) (0.5534,3.6141) (0.5570,3.5906) (0.5607,3.5673) (0.5643,3.5443) (0.5679,3.5216) (0.5716,3.4992) (0.5752,3.4771) (0.5788,3.4553) (0.5825,3.4337) (0.5861,3.4124) (0.5897,3.3914) (0.5934,3.3706) (0.5970,3.3501) (0.6006,3.3299) (0.6043,3.3098) (0.6079,3.2900) (0.6115,3.2705) (0.6152,3.2512) (0.6188,3.2321) (0.6224,3.2132) (0.6261,3.1946) (0.6297,3.1761) (0.6333,3.1579) (0.6370,3.1399) (0.6406,3.1221) (0.6442,3.1045) (0.6479,3.0870) (0.6515,3.0698) (0.6551,3.0528) (0.6588,3.0360) (0.6624,3.0193) (0.6660,3.0028) (0.6697,2.9865) (0.6733,2.9704) (0.6769,2.9545) (0.6806,2.9387) (0.6842,2.9231) (0.6878,2.9076) (0.6915,2.8924) (0.6951,2.8772) (0.6987,2.8623) (0.7024,2.8475) (0.7060,2.8328) (0.7096,2.8183) (0.7133,2.8039) (0.7169,2.7897) (0.7206,2.7757) (0.7242,2.7617) (0.7278,2.7479) (0.7315,2.7343) (0.7351,2.7208) (0.7387,2.7074) (0.7424,2.6941) (0.7460,2.6810) (0.7496,2.6680) (0.7533,2.6551) (0.7569,2.6424) (0.7605,2.6298) (0.7642,2.6173) (0.7678,2.6049) (0.7714,2.5926) (0.7751,2.5804) (0.7787,2.5684) (0.7823,2.5565) (0.7860,2.5446) (0.7896,2.5329) (0.7932,2.5213) (0.7969,2.5098) (0.8005,2.4984) (0.8041,2.4871) (0.8078,2.4760) (0.8114,2.4649) (0.8150,2.4539) (0.8187,2.4430) (0.8223,2.4322) (0.8259,2.4215) (0.8296,2.4109) (0.8332,2.4004) (0.8368,2.3899) (0.8405,2.3796) (0.8441,2.3694) (0.8477,2.3592) (0.8514,2.3491) (0.8550,2.3391) (0.8586,2.3292) (0.8623,2.3194) (0.8659,2.3097) (0.8695,2.3000) (0.8732,2.2905) (0.8768,2.2810) (0.8805,2.2716) (0.8841,2.2622) (0.8877,2.2530) (0.8914,2.2438) (0.8950,2.2347) (0.8986,2.2256) (0.9023,2.2167) (0.9059,2.2078) (0.9095,2.1990) (0.9132,2.1902) (0.9168,2.1815) (0.9204,2.1729) (0.9241,2.1644) (0.9277,2.1559) (0.9313,2.1475) (0.9350,2.1391) (0.9386,2.1308) (0.9422,2.1226) (0.9459,2.1145) (0.9495,2.1064) (0.9531,2.0983) (0.9568,2.0904) (0.9604,2.0825) (0.9640,2.0746) (0.9677,2.0668) (0.9713,2.0591) (0.9749,2.0514) (0.9786,2.0438) (0.9822,2.0362) (0.9858,2.0287) (0.9895,2.0213) (0.9931,2.0139) (0.9967,2.0065) (1.0004,1.9992) (1.0040,1.9920) (1.0076,1.9848) (1.0113,1.9777) (1.0149,1.9706) (1.0185,1.9636) (1.0222,1.9566) (1.0258,1.9497) (1.0294,1.9428) (1.0331,1.9360) (1.0367,1.9292) (1.0404,1.9224) (1.0440,1.9157) (1.0476,1.9091) (1.0513,1.9025) (1.0549,1.8959) (1.0585,1.8894) (1.0622,1.8830) (1.0658,1.8765) (1.0694,1.8702) (1.0731,1.8638) (1.0767,1.8575) (1.0803,1.8513) (1.0840,1.8451) (1.0876,1.8389) (1.0912,1.8328) (1.0949,1.8267) (1.0985,1.8207) (1.1021,1.8147) (1.1058,1.8087) (1.1094,1.8028) (1.1130,1.7969) (1.1167,1.7910) (1.1203,1.7852) (1.1239,1.7795) (1.1276,1.7737) (1.1312,1.7680) (1.1348,1.7624) (1.1385,1.7567) (1.1421,1.7512) (1.1457,1.7456) (1.1494,1.7401) (1.1530,1.7346) (1.1566,1.7291) (1.1603,1.7237) (1.1639,1.7183) (1.1675,1.7130) (1.1712,1.7077) (1.1748,1.7024) (1.1784,1.6972) (1.1821,1.6919) (1.1857,1.6867) (1.1893,1.6816) (1.1930,1.6765) (1.1966,1.6714) (1.2003,1.6663) (1.2039,1.6613) (1.2075,1.6563) (1.2112,1.6513) (1.2148,1.6464) (1.2184,1.6415) (1.2221,1.6366) (1.2257,1.6317) (1.2293,1.6269) (1.2330,1.6221) (1.2366,1.6173) (1.2402,1.6126) (1.2439,1.6079) (1.2475,1.6032) (1.2511,1.5986) (1.2548,1.5939) (1.2584,1.5893) (1.2620,1.5847) (1.2657,1.5802) (1.2693,1.5757) (1.2729,1.5712) (1.2766,1.5667) (1.2802,1.5623) (1.2838,1.5578) (1.2875,1.5534) (1.2911,1.5491) (1.2947,1.5447) (1.2984,1.5404) (1.3020,1.5361) (1.3056,1.5318) (1.3093,1.5276) (1.3129,1.5233) (1.3165,1.5191) (1.3202,1.5150) (1.3238,1.5108) (1.3274,1.5067) (1.3311,1.5025) (1.3347,1.4985) (1.3383,1.4944) (1.3420,1.4903) (1.3456,1.4863) (1.3492,1.4823) (1.3529,1.4783) (1.3565,1.4744) (1.3602,1.4704) (1.3638,1.4665) (1.3674,1.4626) (1.3711,1.4587) (1.3747,1.4549) (1.3783,1.4510) (1.3820,1.4472) (1.3856,1.4434) (1.3892,1.4397) (1.3929,1.4359) (1.3965,1.4322) (1.4001,1.4284) (1.4038,1.4247) (1.4074,1.4211) (1.4110,1.4174) (1.4147,1.4138) (1.4183,1.4101) (1.4219,1.4065) (1.4256,1.4030) (1.4292,1.3994) (1.4328,1.3958) (1.4365,1.3923) (1.4401,1.3888) (1.4437,1.3853) (1.4474,1.3818) (1.4510,1.3784) (1.4546,1.3749) (1.4583,1.3715) (1.4619,1.3681) (1.4655,1.3647) (1.4692,1.3613) (1.4728,1.3580) (1.4764,1.3546) (1.4801,1.3513) (1.4837,1.3480) (1.4873,1.3447) (1.4910,1.3414) (1.4946,1.3381) (1.4982,1.3349) (1.5019,1.3317) (1.5055,1.3285) (1.5091,1.3253) (1.5128,1.3221) (1.5164,1.3189) (1.5201,1.3157) (1.5237,1.3126) (1.5273,1.3095) (1.5310,1.3064) (1.5346,1.3033) (1.5382,1.3002) (1.5419,1.2971) (1.5455,1.2941) (1.5491,1.2911) (1.5528,1.2880) (1.5564,1.2850) (1.5600,1.2820) (1.5637,1.2791) (1.5673,1.2761) (1.5709,1.2731) (1.5746,1.2702) (1.5782,1.2673) (1.5818,1.2644) (1.5855,1.2615) (1.5891,1.2586) (1.5927,1.2557) (1.5964,1.2528) (1.6000,1.2500)};
\addplot[black] coordinates {(0.1609,11.9322) (0.1645,11.6577) (0.1682,11.3950) (0.1718,11.1434) (0.1754,10.9025) (0.1791,10.6712) (0.1827,10.4490) (0.1863,10.2355) (0.1900,10.0303) (0.1936,9.8329) (0.1972,9.6426) (0.2009,9.4591) (0.2045,9.2821) (0.2081,9.1114) (0.2118,8.9468) (0.2154,8.7876) (0.2190,8.6336) (0.2227,8.4845) (0.2263,8.3402) (0.2299,8.2005) (0.2336,8.0654) (0.2372,7.9345) (0.2409,7.8074) (0.2445,7.6840) (0.2481,7.5641) (0.2518,7.4477) (0.2554,7.3345) (0.2590,7.2247) (0.2627,7.1179) (0.2663,7.0143) (0.2699,6.9134) (0.2736,6.8151) (0.2772,6.7193) (0.2808,6.6259) (0.2845,6.5348) (0.2881,6.4460) (0.2917,6.3595) (0.2954,6.2751) (0.2990,6.1928) (0.3026,6.1126) (0.3063,6.0345) (0.3099,5.9582) (0.3135,5.8836) (0.3172,5.8106) (0.3208,5.7392) (0.3244,5.6694) (0.3281,5.6010) (0.3317,5.5341) (0.3353,5.4687) (0.3390,5.4046) (0.3426,5.3420) (0.3462,5.2807) (0.3499,5.2207) (0.3535,5.1620) (0.3571,5.1046) (0.3608,5.0484) (0.3644,4.9935) (0.3680,4.9396) (0.3717,4.8866) (0.3753,4.8347) (0.3789,4.7836) (0.3826,4.7335) (0.3862,4.6843) (0.3898,4.6359) (0.3935,4.5885) (0.3971,4.5419) (0.4008,4.4961) (0.4044,4.4512) (0.4080,4.4070) (0.4117,4.3637) (0.4153,4.3211) (0.4189,4.2793) (0.4226,4.2383) (0.4262,4.1980) (0.4298,4.1584) (0.4335,4.1195) (0.4371,4.0814) (0.4407,4.0439) (0.4444,4.0071) (0.4480,3.9709) (0.4516,3.9352) (0.4553,3.9001) (0.4589,3.8655) (0.4625,3.8313) (0.4662,3.7977) (0.4698,3.7645) (0.4734,3.7319) (0.4771,3.6996) (0.4807,3.6679) (0.4843,3.6366) (0.4880,3.6058) (0.4916,3.5754) (0.4952,3.5454) (0.4989,3.5159) (0.5025,3.4868) (0.5061,3.4582) (0.5098,3.4299) (0.5134,3.4020) (0.5170,3.3746) (0.5207,3.3475) (0.5243,3.3209) (0.5279,3.2946) (0.5316,3.2687) (0.5352,3.2432) (0.5388,3.2180) (0.5425,3.1932) (0.5461,3.1688) (0.5497,3.1447) (0.5534,3.1210) (0.5570,3.0976) (0.5607,3.0745) (0.5643,3.0518) (0.5679,3.0294) (0.5716,3.0073) (0.5752,2.9856) (0.5788,2.9640) (0.5825,2.9427) (0.5861,2.9217) (0.5897,2.9008) (0.5934,2.8802) (0.5970,2.8598) (0.6006,2.8397) (0.6043,2.8197) (0.6079,2.8000) (0.6115,2.7805) (0.6152,2.7612) (0.6188,2.7421) (0.6224,2.7232) (0.6261,2.7045) (0.6297,2.6860) (0.6333,2.6677) (0.6370,2.6497) (0.6406,2.6318) (0.6442,2.6141) (0.6479,2.5966) (0.6515,2.5793) (0.6551,2.5621) (0.6588,2.5452) (0.6624,2.5285) (0.6660,2.5119) (0.6697,2.4955) (0.6733,2.4793) (0.6769,2.4632) (0.6806,2.4474) (0.6842,2.4317) (0.6878,2.4161) (0.6915,2.4008) (0.6951,2.3856) (0.6987,2.3706) (0.7024,2.3557) (0.7060,2.3410) (0.7096,2.3265) (0.7133,2.3121) (0.7169,2.2979) (0.7206,2.2838) (0.7242,2.2699) (0.7278,2.2561) (0.7315,2.2425) (0.7351,2.2290) (0.7387,2.2157) (0.7424,2.2025) (0.7460,2.1895) (0.7496,2.1766) (0.7533,2.1638) (0.7569,2.1512) (0.7605,2.1387) (0.7642,2.1264) (0.7678,2.1141) (0.7714,2.1021) (0.7751,2.0901) (0.7787,2.0783) (0.7823,2.0666) (0.7860,2.0550) (0.7896,2.0436) (0.7932,2.0323) (0.7969,2.0211) (0.8005,2.0100) (0.8041,1.9991) (0.8078,1.9882) (0.8114,1.9774) (0.8150,1.9667) (0.8187,1.9561) (0.8223,1.9456) (0.8259,1.9351) (0.8296,1.9247) (0.8332,1.9144) (0.8368,1.9041) (0.8405,1.8940) (0.8441,1.8839) (0.8477,1.8739) (0.8514,1.8639) (0.8550,1.8540) (0.8586,1.8442) (0.8623,1.8345) (0.8659,1.8248) (0.8695,1.8152) (0.8732,1.8057) (0.8768,1.7963) (0.8805,1.7869) (0.8841,1.7775) (0.8877,1.7683) (0.8914,1.7591) (0.8950,1.7500) (0.8986,1.7409) (0.9023,1.7320) (0.9059,1.7230) (0.9095,1.7142) (0.9132,1.7054) (0.9168,1.6967) (0.9204,1.6880) (0.9241,1.6794) (0.9277,1.6708) (0.9313,1.6624) (0.9350,1.6539) (0.9386,1.6456) (0.9422,1.6373) (0.9459,1.6291) (0.9495,1.6209) (0.9531,1.6128) (0.9568,1.6047) (0.9604,1.5967) (0.9640,1.5888) (0.9677,1.5809) (0.9713,1.5731) (0.9749,1.5653) (0.9786,1.5576) (0.9822,1.5499) (0.9858,1.5423) (0.9895,1.5348) (0.9931,1.5273) (0.9967,1.5198) (1.0004,1.5124) (1.0040,1.5051) (1.0076,1.4978) (1.0113,1.4906) (1.0149,1.4834) (1.0185,1.4763) (1.0222,1.4692) (1.0258,1.4622) (1.0294,1.4552) (1.0331,1.4483) (1.0367,1.4415) (1.0404,1.4346) (1.0440,1.4279) (1.0476,1.4212) (1.0513,1.4145) (1.0549,1.4079) (1.0585,1.4013) (1.0622,1.3948) (1.0658,1.3883) (1.0694,1.3818) (1.0731,1.3755) (1.0767,1.3691) (1.0803,1.3628) (1.0840,1.3566) (1.0876,1.3504) (1.0912,1.3442) (1.0949,1.3381) (1.0985,1.3320) (1.1021,1.3260) (1.1058,1.3200) (1.1094,1.3141) (1.1130,1.3082) (1.1167,1.3023) (1.1203,1.2965) (1.1239,1.2908) (1.1276,1.2850) (1.1312,1.2794) (1.1348,1.2737) (1.1385,1.2681) (1.1421,1.2625) (1.1457,1.2570) (1.1494,1.2515) (1.1530,1.2461) (1.1566,1.2407) (1.1603,1.2353) (1.1639,1.2300) (1.1675,1.2247) (1.1712,1.2195) (1.1748,1.2143) (1.1784,1.2091) (1.1821,1.2040) (1.1857,1.1989) (1.1893,1.1938) (1.1930,1.1888) (1.1966,1.1838) (1.2003,1.1788) (1.2039,1.1739) (1.2075,1.1690) (1.2112,1.1642) (1.2148,1.1594) (1.2184,1.1546) (1.2221,1.1499) (1.2257,1.1452) (1.2293,1.1405) (1.2330,1.1358) (1.2366,1.1312) (1.2402,1.1267) (1.2439,1.1221) (1.2475,1.1176) (1.2511,1.1132) (1.2548,1.1087) (1.2584,1.1043) (1.2620,1.0999) (1.2657,1.0956) (1.2693,1.0913) (1.2729,1.0870) (1.2766,1.0827) (1.2802,1.0785) (1.2838,1.0743) (1.2875,1.0702) (1.2911,1.0660) (1.2947,1.0619) (1.2984,1.0579) (1.3020,1.0538) (1.3056,1.0498) (1.3093,1.0458) (1.3129,1.0419) (1.3165,1.0380) (1.3202,1.0341) (1.3238,1.0302) (1.3274,1.0264) (1.3311,1.0225) (1.3347,1.0188) (1.3383,1.0150) (1.3420,1.0113) (1.3456,1.0076) (1.3492,1.0039) (1.3529,1.0002) (1.3565,0.9966) (1.3602,0.9930) (1.3638,0.9894) (1.3674,0.9858) (1.3711,0.9822) (1.3747,0.9787) (1.3783,0.9751) (1.3820,0.9716) (1.3856,0.9681) (1.3892,0.9646) (1.3929,0.9611) (1.3965,0.9576) (1.4001,0.9541) (1.4038,0.9506) (1.4074,0.9472) (1.4110,0.9437) (1.4147,0.9403) (1.4183,0.9369) (1.4219,0.9335) (1.4256,0.9301) (1.4292,0.9268) (1.4328,0.9234) (1.4365,0.9200) (1.4401,0.9167) (1.4437,0.9134) (1.4474,0.9101) (1.4510,0.9068) (1.4546,0.9035) (1.4583,0.9002) (1.4619,0.8969) (1.4655,0.8937) (1.4692,0.8904) (1.4728,0.8872) (1.4764,0.8840) (1.4801,0.8808) (1.4837,0.8776) (1.4873,0.8744) (1.4910,0.8712) (1.4946,0.8681) (1.4982,0.8649) (1.5019,0.8618) (1.5055,0.8587) (1.5091,0.8555) (1.5128,0.8524) (1.5164,0.8493) (1.5201,0.8463) (1.5237,0.8432) (1.5273,0.8401) (1.5310,0.8371) (1.5346,0.8341) (1.5382,0.8310) (1.5419,0.8280) (1.5455,0.8250) (1.5491,0.8220) (1.5528,0.8190) (1.5564,0.8161) (1.5600,0.8131) (1.5637,0.8102) (1.5673,0.8072) (1.5709,0.8043) (1.5746,0.8014) (1.5782,0.7985) (1.5818,0.7956) (1.5855,0.7927) (1.5891,0.7898) (1.5927,0.7869) (1.5964,0.7841) (1.6000,0.7812)};
\node[anchor=west,font=\small] at (axis cs:1.6,1.550) {$v(T)=2/T$};
\node[anchor=west,font=\small] at (axis cs:1.6,0.831) {$u(T)$};
\node[anchor=south,font=\small,gray!70!black] at (axis cs:1.08,2) {$\varkappa(T)+1$};

\end{axis}
\end{tikzpicture}
\caption{
For sampling   from  $\{0, 1,\ldots, N-1\}$  (respectively, from $\{1,2,\ldots,N\}$), a draw  falling   strictly (respectively, nonstrictly)   below the step curve $\varkappa(T)+1$
 must   be accepted
when  $n\approx TN$ trials  remain.
The lower smooth curve is   the limit stopping value  $u(T)$ for the discrete-uniform problem, and the upper smooth curve
is the limit $v(T)=2/T$ for Moser's problem.
 }\label{Fig1}
\end{figure}
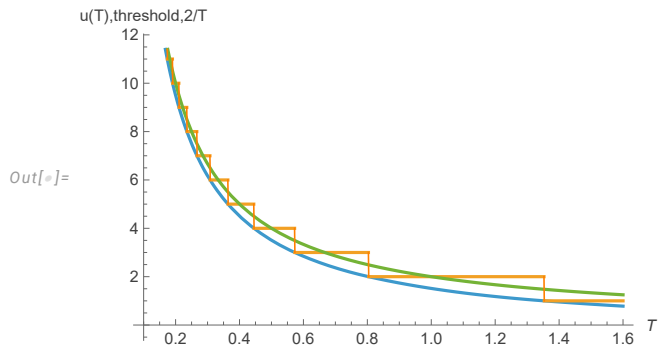

The difference $2/T-u(T)$ between the limit stopping values for the continuous and discrete uniform sampling schemes converges to $1/2$ as $T\to0$, and to $0$ as $T\to\infty$.
The case $T\to\infty$ is straightforward: for large $n/N$,  there will be many zeroes  among the
draws from $\{0,\ldots,N-1\}$,
hence the online strategy is not much different from prophet's.
For small values of $T$ we have the following refinement of (\ref{CompU}).

\begin{prop}\label{P1} Let $n,N\to\infty$ with $n/N\to T$, where $0<T<\sqrt{48}$. Then
\begin{equation}\label{LS}
0\leq \lim_{n/N\to T} \Big(V_{n,N}-\frac{2N}{n}\Big)+\frac{1}{2}=u(T)-\frac{2}{T}+\frac12\leq \frac{T}{24-T^2/2}.
\end{equation}
In particular, the limit difference between the discrete and the continuous problem tends to $-1/2$ as $T\to0$, while $V_{n,N}\,n/(2N)\to Tu(T)/2$, which tends to $1$ as $T\to0$.
\end{prop}

\proof We combine Theorem \ref{T4} with the bounds (\ref{uBounds}) of Lemma \ref{L2}, which apply to $u$ by Theorem \ref{T3}, and note that 
$\frac{2}{T(1-T^2/48)}-\frac{2}{T}=\frac{T}{24-T^2/2}.$
\endpf

The correction $-1/2$ has a transparent meaning: when $T$ is small the accepted draws are large, and the rounding $x\mapsto\lfloor x\rfloor$ of a draw accepted below a high threshold reduces 
the mean loss  by about one half.
Thus the approximation of the discrete problem by Moser's problem holds in the relative sense, $V_{n,N}\sim 2N/n$, when $n/N\to 0$, but not in the sense of absolute differences.
For the original Cayley's problem of sampling without replacement the comparison 
from \cite{Brien} shows that the stopping value does not exceed $V_{n,N}$, so the upper bound $2N/n$ asymptotically applies there as well; we do not pursue the matching lower bound here.

\section{Lindley's minimum rank  problem}

Let $Y_1, Y_2,\ldots$ be independent
random variables, with $Y_j$ having the discrete uniform distribution on $\{1,2,\ldots,j\}$. Lindley's \cite{Lindley} problem with $n$ steps asks
 to find the stopping value
\begin{equation}\label{Lindley}
R_n:=\inf_{\theta\leq n} \me\left(Y_\theta  \frac{n+1}{\theta+1} \right),
\end{equation}
where the infimum is taken over stopping times $\theta$ adapted to the sequence $Y_1,\ldots, Y_n$.
This differs from (\ref{disMP}) in two ways: the range of the distribution changes from trial to trial and there is a time-dependent factor that penalises early stopping.
The first variable has always the smallest possible value, $Y_1=1$, but due to the effect of penalisation  stopping at the first trial yields the same expected loss
$(n+1)/2$ as stopping at the last trial, or
 with any other constant stopping rule $\theta\equiv j$ for $1\leq j\leq n$.

To motivate the objective (\ref{Lindley}) we recall the original Lindley's formulation,
which he introduced as a generalisation of the (classical) `secretary problem' of maximising the probability of the best choice \cite{Ferguson}.
Suppose a choice is to be made from a set  of  $n$ items ranked
$1$ (best) to $n$ (worst).
The items appear one-by-one in random order, and as the $j$th item is inspected  we learn its rank  $Y_j$ {\it relative} to the first $j$ items, but not the overall rank among $n$.
The objective  is to minimise, by way of a stopping rule, the expected overall rank of the accepted item.
By a remarkable combinatorial bijection (known as Lehmer code),  the relative ranks $Y_1,\ldots,Y_n$  with   distribution
as introduced above encode a uniformly random permutation of $n$ integers (overall ranks).
Furthermore,
conditionally on $Y_1,\ldots,Y_{j-1}, Y_j=r$ the expected overall rank of the $j$th item is $r(n+1)/(j+1)$.
Thus  (\ref{Lindley})  is equivalent  to minimising the expected  rank among $n$ items, when the
online observer only learns their relative ranks, as compared to some intrinsic utilities whose distribution might not be known.

Lindley   approached the problem by replacing  the discrete-time optimality equation with    a single differential equation (\ref{odef}) (see \cite[Equation (31)]{Lindley}),
hence obtaining $2$ as an approximate stopping value.
Then   (with reference to a hint  from anonymous  Editor) he proposed the approximate value  $3$, arguing that no stop should be made before  draw $j\approx n/3$, where
the cutoff $1/3$ apparently came from  the balance  equation $2/(1-t)=1/t$.

Chow et al \cite{CMRS}  working directly with the difference equation
succeeded to determine  the limit stopping value
\begin{equation}\label{CMRS}
\lim_{n\to\infty} R_n = \prod_{j=1}^\infty  \left(1+\frac{2}{j}\right)^{\tfrac{1}{j+1}}=3.869\ldots,
\end{equation}
which is a major highlight that stimulated a lot of interest in the discrete-time optimal stopping,
especially on variations  of  the `secretary problem'.
Mucci \cite{Mucci1, Mucci2} advanced the techniques  of differential equations for (\ref{Lindley})  and a wider class  of
problems with rank-dependent  payoffs.   Various  continuous-time   `$n=\infty$'
models  were introduced in \cite{BG, BR, GS1} to extend  the probabilistic interpretation of the limits akin to (\ref{CMRS}).

It goes without saying that the coincidence of   (\ref{edelta})  and (\ref{CMRS})  cannot be accidental,
though  (discrete-uniform)  Moser's and Lindley's problems  appear at a first glance very different. We will link them  using a transformation  of the Poisson limits.

In line with the other parts of this paper, we  represent the relative ranks  via the uniform  random variables as
$Y_{j}=\lfloor j \,U_j\rfloor+1$. Then the arithmetic bounds
  $  j U_j < \lfloor j \,U_j\rfloor+1<  j U_j+1$
endowed with the  penalisation factor become
\begin{equation}\label{ULB}
    j U_j  \,  \frac{n+1}{j+1}< Y_j\, \frac{n+1}{j+1}<  j U_j\,\frac{n+1}{j+1}+\frac{n+1}{j+1}.
\end{equation}
Each of these corresponds to a loss structure for a limit Poisson process.

Consider two point scatters
$${\rm A}_n=\big((j/n, A_j), ~j\in {\mathbb N}),~~~{\rm B}_n=\big((j/n, B_j), ~ j\in {\mathbb N} \big),$$
where the terms
\begin{equation}\label{fr}
A_j:=\big((n+1)U_j\big) \frac{j}{j+1}, ~~ B_j:= jU_j,
\end{equation}
are connected as
\begin{equation}\label{sr}
A_j=    B_j  \frac{n+1}{j+1}.
\end{equation}

\subsection{Lower bound}

Comparing  (\ref{fr}) with scaling  (\ref{SP}) we  note that factor $j/(j+1)$ does not impact the limit, hence we  see that ${\rm A}_n$ is approximable by the homogeneous planar Poisson process
$\Pi$. In the stopping problem for $A_1,\ldots,A_n$ the relevant limit is $\Pi_1$ as we treated in Section \ref{S2.2}. It follows that Lindley's initial approximation  $2$  is an asymptotic  lower bound for $R_n$.

It is much more insightful to view  the bound in terms of  the scatter ${\rm B}_n$ of non-i.i.d. uniforms. Coupled via (\ref{fr}), the variables $B_j$ are smaller than $A_j$, especially for $j/n$ close to zero.
The   measure concentration leads to another, now
inhomogeneous, limit  Poisson process $\widehat{\Pi}$, with the
 rate function $t^{-1 }$ in  ${\mathbb R}_+^2$
(see \cite{GKS} for a different Poisson limit in a  triangular scheme of discrete uniform variables).
The rate function obtains by noting that $j U_j<yj/n$ with probability $y/n$  (for $n>y$).
The process $\widehat{\Pi}$ has infinitely many atoms near every point of the vertical half-axis.
The restriction
 $\widehat{\Pi}_1:=\widehat{\Pi}|_{(0,1]\times{\mathbb R}_+}$ we are interested in still   has this property.
Though the way we label atoms is not really important, we may do this by  restricting further to  times $[\varepsilon, 1]$ (as sufficient for our purpose),
then   label the atoms $(t_i, y_i)$ by increase of $y_i$'s.

In the $\widehat{\Pi}_1$-setting consider the loss structure $(t,y)\mapsto y/t$, where $1/t$ appears as the limit form of penalisation in (\ref{Lindley}), and  also
as the rate.
This twofold role is crucial to identify the stopping problem with that for the homogeneous process $\Pi_1$, so
 the higher  intensity gets compensated by the penalisation.  This somewhat counter-intuitive phenomenon
becomes evident upon noting that  the mapping $(t,y)\mapsto (t, y/t)$ transforms  $\widehat{\Pi}$ into homogeneous $\Pi$.

\subsection{Upper bound}

The optimisation problem associated with the upper bound in  (\ref{ULB}) concerns stopping the Poisson process $\Pi_1$ with loss  structure $(t,x)\mapsto x+1/t$,
which turns interesting by itself.

We denote the stopping value function $h(t)$,  with understanding of the temporal variable $t\in[0,1]$ as the  forward time. By the optimality principle,
with account of
the acceptance condition $x+1/t\leq h(t)$ we obtain
$$
h(t)=dt \int_0^{h(t)-1/t}   \left(x+\frac{1}{t}\right)dx+\left(1-\left(h(t)-\frac{1}{t}\right)dt\right) h(t+dt),
$$
which yields the nonseparable differential equation
$$h'(t)= \frac{h^2(t)}{2}+\frac{1}{2t^2}-\frac{h(t)}{t}.$$
Passing to $g(t)=t h(t)$ reduces this to the separable ODE $2t g'(t)=g^2(t)+1$. Working this out yields the unique solution satisfying    $h(1)=\infty$:
$$
h(t)=\frac{1}{t}\tan\left(\frac{\pi+\log t}{2} \right).
$$
The formula exhibits the correct expansion $h(t)=2/(1-t)+1+o(1-t)$ near $t=1$.

Stopping is not optimal  before  draw value $0$ becomes acceptable at some cutoff time $t_0$.
 The cutoff  is found via continuous fitting from $h(t)=1/t$ as
$t_0=e^{-\pi/2}=0.207\ldots$, therefore
$$h(0)=h(t_0)=e^{\pi/2}=4.810\ldots$$
which is  an asymptotic upper bound for $R_n$ from the right side of (\ref{ULB}).

\subsection{Connection  to stopping of discrete i.i.d. uniform draws}

For Lindley's problem the appropriate limit form is the inhomogeneous Poisson process $\widehat{\Pi}_1$ with loss  structure $(t,y)\mapsto (\lfloor y\rfloor +1)/t$.
In this framework the draws from the box $[k-1,k)$ correspond to the relative rank $k\geq 1$, and they appear according to independent Poisson processes at rate $1/t$.
The limit can be confirmed `box-wise'  as follows. Each value $k$ appears  in position $j\geq k$ of the sequence $Y_1,Y_2,\ldots$ with probability $1/j$,  which
is an inhomogeneous Bernoulli process approximable by the Poisson process with  rate $1/t$.
The piecewise differential equation for the value function in forward time $t\in[0,1]$
was scrutinised  in \cite{Mucci1}, albeit without connecting to a limit model.

To connect with the setting  of Section \ref{S2.2} we pass to a logarithmic time scale by  the time change $t=e^{-T}$,
thus mapping $\widehat{\Pi}$ to the homogeneous $\Pi$.
This works because  the mean number of $\widehat{\Pi}$-atoms falling in  a box  within time interval $e^{-T_2}<t<e^{-T_1}$ is $T_2-T_1$.
That way the limit minimal rank problem is reduced to  optimal stopping of $\Pi_T$, for $T$ not too small (any  $T>\delta_1=1.353\ldots$ would do).
The loss structure becomes   $(t,x)\mapsto (\lfloor x\rfloor +1)e^{T-t}$, where the exponential factor penalises
 early stopping.  Let $w(T)$ be the value function for this problem, that is the infimum of $\me\,(\lfloor\xi_\tau\rfloor+1)e^{T-\tau}$ over stopping times $\tau$ with strict horizon $T$.

\begin{theorem}\label{T5} The stopping value function
$w(T)$   for the Poisson process $\Pi_T$ with the loss structure $(t,x)\mapsto (\lfloor x\rfloor +1)e^{T-t}$
is
\begin{equation}\label{wu}
w(T)=e^T u(T),
\end{equation}
for $u(T)$ as in {\rm Theorem \ref{T3}}.
The optimal stopping time in this problem is
$$
\tilde{\nu}(T):=\min\{t_i: t_i\geq (T-\delta_k)_+~~{\rm for~~}k=\lfloor x_i\rfloor+1\},
$$
with cutoffs {\rm (\ref{delta-cut})}.
\end{theorem}

\proof Let $\Pi^{(0)}$ be a unit-rate Poisson process in $[0,T]\times[0,1)$, independent of $\Pi$, and let $\tau_0$ be the time of its first atom ($\tau_0:=\infty$ if there is none), so that $\prob(\tau_0>t)=e^{-t}$. The superposition
$$\Pi^\uparrow:=\Pi^{(0)}\cup\{(t_i,x_i+1)\}$$
of $\Pi^{(0)}$ and the lifted atoms of $\Pi_T$ is again a unit-rate Poisson process in $[0,T]\times{\mathbb R}_+$, hence $u(T)$ is the stopping value of the problem (\ref{dePo}) for $\Pi^\uparrow$. In this problem the atoms of $\Pi^{(0)}$ incur zero loss, so replacing a stopping time $\rho$ by $\min(\rho,\tau_0)$ does not increase the loss. Before $\tau_0$ the process $\Pi^\uparrow$ conveys the same information as $\Pi$, therefore $\min(\rho,\tau_0)=\min(\tau,\tau_0)$ for some stopping time $\tau$ with strict horizon $T$ adapted to $\Pi$ (namely, $\tau$ acts as $\rho$ does when $\Pi^{(0)}$ has no atoms). Conversely, each such $\tau$ defines in this way a stopping time for $\Pi^\uparrow$. The rule $\min(\tau,\tau_0)$ incurs the loss $(\lfloor\xi_\tau\rfloor+1)1(\tau<\tau_0)$, and by independence of $\tau_0$ and $\Pi$
$$
\me\big[(\lfloor\xi_\tau\rfloor+1)1(\tau<\tau_0)\big]=\me\big[(\lfloor\xi_\tau\rfloor+1)e^{-\tau}\big]=e^{-T}\,\me\big[(\lfloor\xi_\tau\rfloor+1)e^{T-\tau}\big].
$$
Taking the infimum over $\tau$ yields $u(T)=e^{-T}w(T)$.

The optimal rule $\nu(T)$ for $\Pi^\uparrow$ accepts every atom of $\Pi^{(0)}$, and it accepts a lifted atom $(t_i,x_i+1)$ with $\lfloor x_i\rfloor+1=k$ if and only if $t_i\geq(T-\delta_k)_+$. Thus $\nu(T)=\min(\tilde\nu(T),\tau_0)$, and the optimality of $\tilde\nu(T)$ follows from the above correspondence.
\endpf

In particular, the two optimal rules are related through the coupling $\Pi^\uparrow$ rather than through acting on the same realisation of $\Pi$: strategy $\tilde\nu(T)$ on $\Pi$ behaves as $\nu(T)$ on the process obtained by lifting all atoms by one unit, with the additional zero-loss atoms of the bottom box replaced by the exponential discount $e^{-\tau}$. Acting on the same realisation of $\Pi$ the two rules may stop at different atoms; for instance, if $\delta_2<T<\delta_1$ then before time $T-\delta_2$ the rule $\nu(T)$ accepts values in $[0,2)$ while $\tilde\nu(T)$ accepts only values in $[0,1)$.

As a consistency check, $w=e^Tu$ satisfies the piecewise differential equation
$$w'(T)=-\sum_{j=1}^\infty(w(T)-e^Tj)_+, ~~~w(0)=\infty,$$
which is the optimality equation for the loss structure $(t,x)\mapsto(\lfloor x\rfloor+1)e^{T-t}$; it involves the penalising factor and has no $j=0$ term. In the range $e^Tk<w(T)<e^T(k+1)$, that is $\delta_{k+1}<T<\delta_k$, the solution is $w(T)=\frac{k}{2}e^T+e^{-kT}C$, and fitting the boundary conditions recovers the recursion (\ref{reccut}).

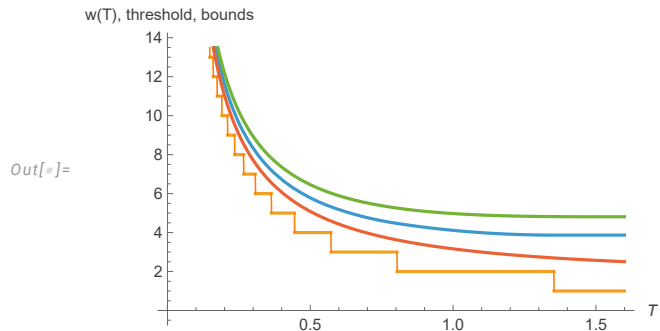
\begin{figure}[htbp]
\centering
\begin{tikzpicture}
\begin{axis}[width=0.8\textwidth,height=0.5\textwidth,xmin=0.15,xmax=1.6,ymin=0,ymax=14,
 axis lines=left, xlabel={$T$}, xlabel style={at={(axis description cs:1.02,0)},anchor=west},
 ytick={0,2,...,14}, xtick={0.2,0.4,...,1.6}, tick label style={font=\small},
 every axis plot/.append style={line width=0.9pt}, clip=false]
\addplot[gray!70!black, line width=0.7pt, const plot, mark=none] coordinates {(0.1500,13) (0.1600,12) (0.1739,11) (0.1905,10) (0.2106,9) (0.2354,8) (0.2668,7) (0.3079,6) (0.3640,5) (0.4451,4) (0.5728,3) (0.8038,2) (1.3531,1) (1.6000,1)};
\addplot[black, densely dotted, line width=1.1pt] coordinates {(0.1718,13.7892) (0.1754,13.5513) (0.1791,13.3232) (0.1827,13.1043) (0.1863,12.8941) (0.1900,12.6920) (0.1936,12.4976) (0.1972,12.3106) (0.2009,12.1304) (0.2045,11.9568) (0.2081,11.7893) (0.2118,11.6277) (0.2154,11.4717) (0.2190,11.3209) (0.2227,11.1752) (0.2263,11.0343) (0.2299,10.8979) (0.2336,10.7659) (0.2372,10.6380) (0.2409,10.5141) (0.2445,10.3940) (0.2481,10.2775) (0.2518,10.1645) (0.2554,10.0548) (0.2590,9.9483) (0.2627,9.8448) (0.2663,9.7442) (0.2699,9.6464) (0.2736,9.5513) (0.2772,9.4588) (0.2808,9.3688) (0.2845,9.2812) (0.2881,9.1959) (0.2917,9.1128) (0.2954,9.0318) (0.2990,8.9529) (0.3026,8.8759) (0.3063,8.8009) (0.3099,8.7277) (0.3135,8.6563) (0.3172,8.5866) (0.3208,8.5186) (0.3244,8.4522) (0.3281,8.3873) (0.3317,8.3239) (0.3353,8.2620) (0.3390,8.2015) (0.3426,8.1423) (0.3462,8.0845) (0.3499,8.0279) (0.3535,7.9726) (0.3571,7.9185) (0.3608,7.8655) (0.3644,7.8137) (0.3680,7.7630) (0.3717,7.7133) (0.3753,7.6647) (0.3789,7.6170) (0.3826,7.5704) (0.3862,7.5246) (0.3898,7.4798) (0.3935,7.4359) (0.3971,7.3929) (0.4008,7.3507) (0.4044,7.3094) (0.4080,7.2688) (0.4117,7.2291) (0.4153,7.1901) (0.4189,7.1518) (0.4226,7.1142) (0.4262,7.0774) (0.4298,7.0413) (0.4335,7.0058) (0.4371,6.9709) (0.4407,6.9367) (0.4444,6.9032) (0.4480,6.8702) (0.4516,6.8378) (0.4553,6.8060) (0.4589,6.7748) (0.4625,6.7441) (0.4662,6.7139) (0.4698,6.6843) (0.4734,6.6552) (0.4771,6.6266) (0.4807,6.5985) (0.4843,6.5708) (0.4880,6.5437) (0.4916,6.5169) (0.4952,6.4907) (0.4989,6.4648) (0.5025,6.4394) (0.5061,6.4145) (0.5098,6.3899) (0.5134,6.3657) (0.5170,6.3420) (0.5207,6.3186) (0.5243,6.2956) (0.5279,6.2729) (0.5316,6.2507) (0.5352,6.2287) (0.5388,6.2072) (0.5425,6.1859) (0.5461,6.1650) (0.5497,6.1445) (0.5534,6.1242) (0.5570,6.1043) (0.5607,6.0846) (0.5643,6.0653) (0.5679,6.0463) (0.5716,6.0276) (0.5752,6.0091) (0.5788,5.9910) (0.5825,5.9731) (0.5861,5.9554) (0.5897,5.9381) (0.5934,5.9210) (0.5970,5.9042) (0.6006,5.8876) (0.6043,5.8712) (0.6079,5.8551) (0.6115,5.8393) (0.6152,5.8237) (0.6188,5.8083) (0.6224,5.7931) (0.6261,5.7781) (0.6297,5.7634) (0.6333,5.7489) (0.6370,5.7346) (0.6406,5.7205) (0.6442,5.7066) (0.6479,5.6929) (0.6515,5.6794) (0.6551,5.6661) (0.6588,5.6530) (0.6624,5.6401) (0.6660,5.6273) (0.6697,5.6148) (0.6733,5.6024) (0.6769,5.5902) (0.6806,5.5782) (0.6842,5.5663) (0.6878,5.5546) (0.6915,5.5431) (0.6951,5.5317) (0.6987,5.5205) (0.7024,5.5095) (0.7060,5.4986) (0.7096,5.4878) (0.7133,5.4773) (0.7169,5.4668) (0.7206,5.4565) (0.7242,5.4464) (0.7278,5.4364) (0.7315,5.4265) (0.7351,5.4167) (0.7387,5.4071) (0.7424,5.3977) (0.7460,5.3883) (0.7496,5.3791) (0.7533,5.3701) (0.7569,5.3611) (0.7605,5.3523) (0.7642,5.3436) (0.7678,5.3350) (0.7714,5.3265) (0.7751,5.3182) (0.7787,5.3099) (0.7823,5.3018) (0.7860,5.2938) (0.7896,5.2859) (0.7932,5.2781) (0.7969,5.2704) (0.8005,5.2629) (0.8041,5.2554) (0.8078,5.2480) (0.8114,5.2408) (0.8150,5.2336) (0.8187,5.2265) (0.8223,5.2196) (0.8259,5.2127) (0.8296,5.2059) (0.8332,5.1992) (0.8368,5.1926) (0.8405,5.1861) (0.8441,5.1797) (0.8477,5.1734) (0.8514,5.1672) (0.8550,5.1610) (0.8586,5.1550) (0.8623,5.1490) (0.8659,5.1431) (0.8695,5.1373) (0.8732,5.1316) (0.8768,5.1259) (0.8805,5.1204) (0.8841,5.1149) (0.8877,5.1095) (0.8914,5.1042) (0.8950,5.0989) (0.8986,5.0937) (0.9023,5.0886) (0.9059,5.0836) (0.9095,5.0786) (0.9132,5.0737) (0.9168,5.0689) (0.9204,5.0641) (0.9241,5.0594) (0.9277,5.0548) (0.9313,5.0503) (0.9350,5.0458) (0.9386,5.0414) (0.9422,5.0370) (0.9459,5.0327) (0.9495,5.0285) (0.9531,5.0243) (0.9568,5.0202) (0.9604,5.0162) (0.9640,5.0122) (0.9677,5.0082) (0.9713,5.0044) (0.9749,5.0006) (0.9786,4.9968) (0.9822,4.9931) (0.9858,4.9895) (0.9895,4.9859) (0.9931,4.9823) (0.9967,4.9789) (1.0004,4.9754) (1.0040,4.9721) (1.0076,4.9687) (1.0113,4.9655) (1.0149,4.9622) (1.0185,4.9591) (1.0222,4.9559) (1.0258,4.9529) (1.0294,4.9498) (1.0331,4.9469) (1.0367,4.9439) (1.0404,4.9411) (1.0440,4.9382) (1.0476,4.9354) (1.0513,4.9327) (1.0549,4.9300) (1.0585,4.9273) (1.0622,4.9247) (1.0658,4.9221) (1.0694,4.9196) (1.0731,4.9171) (1.0767,4.9147) (1.0803,4.9123) (1.0840,4.9099) (1.0876,4.9076) (1.0912,4.9053) (1.0949,4.9030) (1.0985,4.9008) (1.1021,4.8987) (1.1058,4.8965) (1.1094,4.8944) (1.1130,4.8924) (1.1167,4.8903) (1.1203,4.8884) (1.1239,4.8864) (1.1276,4.8845) (1.1312,4.8826) (1.1348,4.8808) (1.1385,4.8789) (1.1421,4.8772) (1.1457,4.8754) (1.1494,4.8737) (1.1530,4.8720) (1.1566,4.8704) (1.1603,4.8688) (1.1639,4.8672) (1.1675,4.8656) (1.1712,4.8641) (1.1748,4.8626) (1.1784,4.8611) (1.1821,4.8597) (1.1857,4.8583) (1.1893,4.8569) (1.1930,4.8555) (1.1966,4.8542) (1.2003,4.8529) (1.2039,4.8516) (1.2075,4.8504) (1.2112,4.8492) (1.2148,4.8480) (1.2184,4.8468) (1.2221,4.8457) (1.2257,4.8446) (1.2293,4.8435) (1.2330,4.8424) (1.2366,4.8414) (1.2402,4.8403) (1.2439,4.8394) (1.2475,4.8384) (1.2511,4.8374) (1.2548,4.8365) (1.2584,4.8356) (1.2620,4.8347) (1.2657,4.8339) (1.2693,4.8330) (1.2729,4.8322) (1.2766,4.8314) (1.2802,4.8306) (1.2838,4.8299) (1.2875,4.8291) (1.2911,4.8284) (1.2947,4.8277) (1.2984,4.8270) (1.3020,4.8264) (1.3056,4.8257) (1.3093,4.8251) (1.3129,4.8245) (1.3165,4.8239) (1.3202,4.8233) (1.3238,4.8228) (1.3274,4.8222) (1.3311,4.8217) (1.3347,4.8212) (1.3383,4.8207) (1.3420,4.8202) (1.3456,4.8198) (1.3492,4.8193) (1.3529,4.8189) (1.3565,4.8185) (1.3602,4.8181) (1.3638,4.8177) (1.3674,4.8173) (1.3711,4.8169) (1.3747,4.8166) (1.3783,4.8163) (1.3820,4.8159) (1.3856,4.8156) (1.3892,4.8153) (1.3929,4.8150) (1.3965,4.8148) (1.4001,4.8145) (1.4038,4.8142) (1.4074,4.8140) (1.4110,4.8138) (1.4147,4.8135) (1.4183,4.8133) (1.4219,4.8131) (1.4256,4.8129) (1.4292,4.8128) (1.4328,4.8126) (1.4365,4.8124) (1.4401,4.8123) (1.4437,4.8121) (1.4474,4.8120) (1.4510,4.8119) (1.4546,4.8117) (1.4583,4.8116) (1.4619,4.8115) (1.4655,4.8114) (1.4692,4.8113) (1.4728,4.8112) (1.4764,4.8112) (1.4801,4.8111) (1.4837,4.8110) (1.4873,4.8109) (1.4910,4.8109) (1.4946,4.8108) (1.4982,4.8108) (1.5019,4.8107) (1.5055,4.8107) (1.5091,4.8107) (1.5128,4.8106) (1.5164,4.8106) (1.5201,4.8106) (1.5237,4.8106) (1.5273,4.8105) (1.5310,4.8105) (1.5346,4.8105) (1.5382,4.8105) (1.5419,4.8105) (1.5455,4.8105) (1.5491,4.8105) (1.5528,4.8105) (1.5564,4.8105) (1.5600,4.8105) (1.5637,4.8105) (1.5673,4.8105) (1.5709,4.8105) (1.5746,4.8105) (1.5782,4.8105) (1.5818,4.8105) (1.5855,4.8105) (1.5891,4.8105) (1.5927,4.8105) (1.5964,4.8105) (1.6000,4.8105)};
\addplot[black] coordinates {(0.1645,13.7427) (0.1682,13.4818) (0.1718,13.2322) (0.1754,12.9933) (0.1791,12.7639) (0.1827,12.5436) (0.1863,12.3320) (0.1900,12.1288) (0.1936,11.9335) (0.1972,11.7451) (0.2009,11.5635) (0.2045,11.3884) (0.2081,11.2197) (0.2118,11.0572) (0.2154,10.9000) (0.2190,10.7479) (0.2227,10.6007) (0.2263,10.4584) (0.2299,10.3207) (0.2336,10.1875) (0.2372,10.0587) (0.2409,9.9337) (0.2445,9.8122) (0.2481,9.6943) (0.2518,9.5798) (0.2554,9.4686) (0.2590,9.3607) (0.2627,9.2560) (0.2663,9.1544) (0.2699,9.0557) (0.2736,8.9594) (0.2772,8.8656) (0.2808,8.7742) (0.2845,8.6851) (0.2881,8.5983) (0.2917,8.5137) (0.2954,8.4313) (0.2990,8.3511) (0.3026,8.2729) (0.3063,8.1969) (0.3099,8.1227) (0.3135,8.0503) (0.3172,7.9793) (0.3208,7.9100) (0.3244,7.8421) (0.3281,7.7758) (0.3317,7.7109) (0.3353,7.6475) (0.3390,7.5854) (0.3426,7.5248) (0.3462,7.4655) (0.3499,7.4076) (0.3535,7.3510) (0.3571,7.2957) (0.3608,7.2417) (0.3644,7.1889) (0.3680,7.1372) (0.3717,7.0865) (0.3753,7.0366) (0.3789,6.9877) (0.3826,6.9396) (0.3862,6.8925) (0.3898,6.8462) (0.3935,6.8008) (0.3971,6.7562) (0.4008,6.7124) (0.4044,6.6695) (0.4080,6.6274) (0.4117,6.5862) (0.4153,6.5457) (0.4189,6.5060) (0.4226,6.4670) (0.4262,6.4288) (0.4298,6.3914) (0.4335,6.3547) (0.4371,6.3188) (0.4407,6.2835) (0.4444,6.2490) (0.4480,6.2151) (0.4516,6.1818) (0.4553,6.1489) (0.4589,6.1164) (0.4625,6.0845) (0.4662,6.0530) (0.4698,6.0220) (0.4734,5.9915) (0.4771,5.9614) (0.4807,5.9318) (0.4843,5.9026) (0.4880,5.8738) (0.4916,5.8455) (0.4952,5.8177) (0.4989,5.7902) (0.5025,5.7632) (0.5061,5.7367) (0.5098,5.7105) (0.5134,5.6847) (0.5170,5.6594) (0.5207,5.6345) (0.5243,5.6100) (0.5279,5.5858) (0.5316,5.5621) (0.5352,5.5387) (0.5388,5.5158) (0.5425,5.4932) (0.5461,5.4710) (0.5497,5.4492) (0.5534,5.4278) (0.5570,5.4067) (0.5607,5.3860) (0.5643,5.3657) (0.5679,5.3457) (0.5716,5.3261) (0.5752,5.3068) (0.5788,5.2877) (0.5825,5.2688) (0.5861,5.2501) (0.5897,5.2317) (0.5934,5.2134) (0.5970,5.1953) (0.6006,5.1775) (0.6043,5.1598) (0.6079,5.1424) (0.6115,5.1251) (0.6152,5.1081) (0.6188,5.0912) (0.6224,5.0745) (0.6261,5.0581) (0.6297,5.0418) (0.6333,5.0257) (0.6370,5.0098) (0.6406,4.9941) (0.6442,4.9786) (0.6479,4.9633) (0.6515,4.9481) (0.6551,4.9332) (0.6588,4.9184) (0.6624,4.9038) (0.6660,4.8894) (0.6697,4.8752) (0.6733,4.8611) (0.6769,4.8473) (0.6806,4.8336) (0.6842,4.8201) (0.6878,4.8067) (0.6915,4.7936) (0.6951,4.7806) (0.6987,4.7678) (0.7024,4.7551) (0.7060,4.7427) (0.7096,4.7304) (0.7133,4.7182) (0.7169,4.7063) (0.7206,4.6945) (0.7242,4.6829) (0.7278,4.6714) (0.7315,4.6601) (0.7351,4.6490) (0.7387,4.6380) (0.7424,4.6272) (0.7460,4.6165) (0.7496,4.6061) (0.7533,4.5957) (0.7569,4.5856) (0.7605,4.5756) (0.7642,4.5657) (0.7678,4.5560) (0.7714,4.5465) (0.7751,4.5371) (0.7787,4.5278) (0.7823,4.5188) (0.7860,4.5098) (0.7896,4.5011) (0.7932,4.4924) (0.7969,4.4839) (0.8005,4.4756) (0.8041,4.4674) (0.8078,4.4594) (0.8114,4.4513) (0.8150,4.4434) (0.8187,4.4355) (0.8223,4.4276) (0.8259,4.4198) (0.8296,4.4121) (0.8332,4.4044) (0.8368,4.3968) (0.8405,4.3892) (0.8441,4.3817) (0.8477,4.3743) (0.8514,4.3669) (0.8550,4.3596) (0.8586,4.3523) (0.8623,4.3451) (0.8659,4.3379) (0.8695,4.3309) (0.8732,4.3238) (0.8768,4.3168) (0.8805,4.3099) (0.8841,4.3030) (0.8877,4.2962) (0.8914,4.2895) (0.8950,4.2828) (0.8986,4.2761) (0.9023,4.2695) (0.9059,4.2630) (0.9095,4.2565) (0.9132,4.2501) (0.9168,4.2438) (0.9204,4.2375) (0.9241,4.2312) (0.9277,4.2250) (0.9313,4.2189) (0.9350,4.2128) (0.9386,4.2068) (0.9422,4.2008) (0.9459,4.1949) (0.9495,4.1890) (0.9531,4.1832) (0.9568,4.1775) (0.9604,4.1718) (0.9640,4.1661) (0.9677,4.1606) (0.9713,4.1550) (0.9749,4.1496) (0.9786,4.1441) (0.9822,4.1388) (0.9858,4.1335) (0.9895,4.1282) (0.9931,4.1230) (0.9967,4.1179) (1.0004,4.1128) (1.0040,4.1077) (1.0076,4.1028) (1.0113,4.0978) (1.0149,4.0930) (1.0185,4.0881) (1.0222,4.0834) (1.0258,4.0787) (1.0294,4.0740) (1.0331,4.0694) (1.0367,4.0648) (1.0404,4.0603) (1.0440,4.0559) (1.0476,4.0515) (1.0513,4.0472) (1.0549,4.0429) (1.0585,4.0387) (1.0622,4.0345) (1.0658,4.0304) (1.0694,4.0263) (1.0731,4.0223) (1.0767,4.0183) (1.0803,4.0144) (1.0840,4.0105) (1.0876,4.0067) (1.0912,4.0030) (1.0949,3.9993) (1.0985,3.9956) (1.1021,3.9921) (1.1058,3.9885) (1.1094,3.9850) (1.1130,3.9816) (1.1167,3.9782) (1.1203,3.9749) (1.1239,3.9716) (1.1276,3.9684) (1.1312,3.9652) (1.1348,3.9621) (1.1385,3.9590) (1.1421,3.9560) (1.1457,3.9530) (1.1494,3.9501) (1.1530,3.9473) (1.1566,3.9445) (1.1603,3.9417) (1.1639,3.9390) (1.1675,3.9363) (1.1712,3.9337) (1.1748,3.9312) (1.1784,3.9287) (1.1821,3.9263) (1.1857,3.9239) (1.1893,3.9215) (1.1930,3.9192) (1.1966,3.9170) (1.2003,3.9148) (1.2039,3.9127) (1.2075,3.9106) (1.2112,3.9086) (1.2148,3.9066) (1.2184,3.9047) (1.2221,3.9028) (1.2257,3.9010) (1.2293,3.8992) (1.2330,3.8975) (1.2366,3.8958) (1.2402,3.8942) (1.2439,3.8926) (1.2475,3.8911) (1.2511,3.8897) (1.2548,3.8883) (1.2584,3.8869) (1.2620,3.8856) (1.2657,3.8843) (1.2693,3.8831) (1.2729,3.8820) (1.2766,3.8809) (1.2802,3.8798) (1.2838,3.8788) (1.2875,3.8779) (1.2911,3.8770) (1.2947,3.8761) (1.2984,3.8753) (1.3020,3.8746) (1.3056,3.8739) (1.3093,3.8732) (1.3129,3.8726) (1.3165,3.8721) (1.3202,3.8716) (1.3238,3.8712) (1.3274,3.8708) (1.3311,3.8705) (1.3347,3.8702) (1.3383,3.8699) (1.3420,3.8698) (1.3456,3.8696) (1.3492,3.8695) (1.3529,3.8695) (1.3565,3.8695) (1.3602,3.8695) (1.3638,3.8695) (1.3674,3.8695) (1.3711,3.8695) (1.3747,3.8695) (1.3783,3.8695) (1.3820,3.8695) (1.3856,3.8695) (1.3892,3.8695) (1.3929,3.8695) (1.3965,3.8695) (1.4001,3.8695) (1.4038,3.8695) (1.4074,3.8695) (1.4110,3.8695) (1.4147,3.8695) (1.4183,3.8695) (1.4219,3.8695) (1.4256,3.8695) (1.4292,3.8695) (1.4328,3.8695) (1.4365,3.8695) (1.4401,3.8695) (1.4437,3.8695) (1.4474,3.8695) (1.4510,3.8695) (1.4546,3.8695) (1.4583,3.8695) (1.4619,3.8695) (1.4655,3.8695) (1.4692,3.8695) (1.4728,3.8695) (1.4764,3.8695) (1.4801,3.8695) (1.4837,3.8695) (1.4873,3.8695) (1.4910,3.8695) (1.4946,3.8695) (1.4982,3.8695) (1.5019,3.8695) (1.5055,3.8695) (1.5091,3.8695) (1.5128,3.8695) (1.5164,3.8695) (1.5201,3.8695) (1.5237,3.8695) (1.5273,3.8695) (1.5310,3.8695) (1.5346,3.8695) (1.5382,3.8695) (1.5419,3.8695) (1.5455,3.8695) (1.5491,3.8695) (1.5528,3.8695) (1.5564,3.8695) (1.5600,3.8695) (1.5637,3.8695) (1.5673,3.8695) (1.5709,3.8695) (1.5746,3.8695) (1.5782,3.8695) (1.5818,3.8695) (1.5855,3.8695) (1.5891,3.8695) (1.5927,3.8695) (1.5964,3.8695) (1.6000,3.8695)};
\addplot[black, dashed] coordinates {(0.1573,13.7433) (0.1609,13.4567) (0.1645,13.1828) (0.1682,12.9207) (0.1718,12.6698) (0.1754,12.4292) (0.1791,12.1985) (0.1827,11.9769) (0.1863,11.7641) (0.1900,11.5593) (0.1936,11.3623) (0.1972,11.1726) (0.2009,10.9898) (0.2045,10.8135) (0.2081,10.6433) (0.2118,10.4791) (0.2154,10.3203) (0.2190,10.1669) (0.2227,10.0185) (0.2263,9.8749) (0.2299,9.7358) (0.2336,9.6011) (0.2372,9.4706) (0.2409,9.3440) (0.2445,9.2211) (0.2481,9.1019) (0.2518,8.9862) (0.2554,8.8737) (0.2590,8.7645) (0.2627,8.6582) (0.2663,8.5549) (0.2699,8.4544) (0.2736,8.3566) (0.2772,8.2613) (0.2808,8.1686) (0.2845,8.0782) (0.2881,7.9901) (0.2917,7.9042) (0.2954,7.8205) (0.2990,7.7388) (0.3026,7.6591) (0.3063,7.5812) (0.3099,7.5053) (0.3135,7.4311) (0.3172,7.3586) (0.3208,7.2878) (0.3244,7.2185) (0.3281,7.1508) (0.3317,7.0847) (0.3353,7.0199) (0.3390,6.9566) (0.3426,6.8946) (0.3462,6.8339) (0.3499,6.7745) (0.3535,6.7164) (0.3571,6.6594) (0.3608,6.6036) (0.3644,6.5489) (0.3680,6.4953) (0.3717,6.4428) (0.3753,6.3913) (0.3789,6.3408) (0.3826,6.2913) (0.3862,6.2427) (0.3898,6.1950) (0.3935,6.1482) (0.3971,6.1023) (0.4008,6.0572) (0.4044,6.0130) (0.4080,5.9695) (0.4117,5.9269) (0.4153,5.8849) (0.4189,5.8438) (0.4226,5.8033) (0.4262,5.7636) (0.4298,5.7245) (0.4335,5.6861) (0.4371,5.6483) (0.4407,5.6112) (0.4444,5.5747) (0.4480,5.5388) (0.4516,5.5034) (0.4553,5.4687) (0.4589,5.4345) (0.4625,5.4008) (0.4662,5.3677) (0.4698,5.3351) (0.4734,5.3031) (0.4771,5.2715) (0.4807,5.2404) (0.4843,5.2098) (0.4880,5.1796) (0.4916,5.1499) (0.4952,5.1207) (0.4989,5.0918) (0.5025,5.0635) (0.5061,5.0355) (0.5098,5.0079) (0.5134,4.9807) (0.5170,4.9539) (0.5207,4.9275) (0.5243,4.9015) (0.5279,4.8759) (0.5316,4.8506) (0.5352,4.8256) (0.5388,4.8010) (0.5425,4.7767) (0.5461,4.7528) (0.5497,4.7292) (0.5534,4.7059) (0.5570,4.6829) (0.5607,4.6602) (0.5643,4.6379) (0.5679,4.6158) (0.5716,4.5940) (0.5752,4.5725) (0.5788,4.5512) (0.5825,4.5303) (0.5861,4.5096) (0.5897,4.4891) (0.5934,4.4690) (0.5970,4.4490) (0.6006,4.4294) (0.6043,4.4099) (0.6079,4.3907) (0.6115,4.3718) (0.6152,4.3531) (0.6188,4.3346) (0.6224,4.3163) (0.6261,4.2982) (0.6297,4.2804) (0.6333,4.2628) (0.6370,4.2453) (0.6406,4.2281) (0.6442,4.2111) (0.6479,4.1943) (0.6515,4.1776) (0.6551,4.1612) (0.6588,4.1450) (0.6624,4.1289) (0.6660,4.1130) (0.6697,4.0973) (0.6733,4.0818) (0.6769,4.0664) (0.6806,4.0512) (0.6842,4.0362) (0.6878,4.0214) (0.6915,4.0067) (0.6951,3.9922) (0.6987,3.9778) (0.7024,3.9636) (0.7060,3.9495) (0.7096,3.9356) (0.7133,3.9218) (0.7169,3.9082) (0.7206,3.8947) (0.7242,3.8814) (0.7278,3.8682) (0.7315,3.8551) (0.7351,3.8422) (0.7387,3.8294) (0.7424,3.8167) (0.7460,3.8042) (0.7496,3.7918) (0.7533,3.7795) (0.7569,3.7673) (0.7605,3.7553) (0.7642,3.7434) (0.7678,3.7316) (0.7714,3.7199) (0.7751,3.7083) (0.7787,3.6969) (0.7823,3.6855) (0.7860,3.6743) (0.7896,3.6632) (0.7932,3.6522) (0.7969,3.6413) (0.8005,3.6304) (0.8041,3.6197) (0.8078,3.6091) (0.8114,3.5986) (0.8150,3.5882) (0.8187,3.5779) (0.8223,3.5677) (0.8259,3.5576) (0.8296,3.5476) (0.8332,3.5376) (0.8368,3.5278) (0.8405,3.5181) (0.8441,3.5084) (0.8477,3.4988) (0.8514,3.4893) (0.8550,3.4799) (0.8586,3.4706) (0.8623,3.4614) (0.8659,3.4522) (0.8695,3.4432) (0.8732,3.4342) (0.8768,3.4253) (0.8805,3.4164) (0.8841,3.4077) (0.8877,3.3990) (0.8914,3.3904) (0.8950,3.3819) (0.8986,3.3734) (0.9023,3.3650) (0.9059,3.3567) (0.9095,3.3485) (0.9132,3.3403) (0.9168,3.3322) (0.9204,3.3242) (0.9241,3.3162) (0.9277,3.3083) (0.9313,3.3005) (0.9350,3.2927) (0.9386,3.2850) (0.9422,3.2774) (0.9459,3.2698) (0.9495,3.2623) (0.9531,3.2548) (0.9568,3.2475) (0.9604,3.2401) (0.9640,3.2329) (0.9677,3.2256) (0.9713,3.2185) (0.9749,3.2114) (0.9786,3.2043) (0.9822,3.1974) (0.9858,3.1904) (0.9895,3.1836) (0.9931,3.1767) (0.9967,3.1700) (1.0004,3.1633) (1.0040,3.1566) (1.0076,3.1500) (1.0113,3.1434) (1.0149,3.1369) (1.0185,3.1305) (1.0222,3.1241) (1.0258,3.1177) (1.0294,3.1114) (1.0331,3.1051) (1.0367,3.0989) (1.0404,3.0928) (1.0440,3.0867) (1.0476,3.0806) (1.0513,3.0746) (1.0549,3.0686) (1.0585,3.0626) (1.0622,3.0567) (1.0658,3.0509) (1.0694,3.0451) (1.0731,3.0393) (1.0767,3.0336) (1.0803,3.0279) (1.0840,3.0223) (1.0876,3.0167) (1.0912,3.0112) (1.0949,3.0056) (1.0985,3.0002) (1.1021,2.9947) (1.1058,2.9893) (1.1094,2.9840) (1.1130,2.9787) (1.1167,2.9734) (1.1203,2.9682) (1.1239,2.9630) (1.1276,2.9578) (1.1312,2.9527) (1.1348,2.9476) (1.1385,2.9425) (1.1421,2.9375) (1.1457,2.9325) (1.1494,2.9276) (1.1530,2.9226) (1.1566,2.9178) (1.1603,2.9129) (1.1639,2.9081) (1.1675,2.9033) (1.1712,2.8986) (1.1748,2.8938) (1.1784,2.8892) (1.1821,2.8845) (1.1857,2.8799) (1.1893,2.8753) (1.1930,2.8707) (1.1966,2.8662) (1.2003,2.8617) (1.2039,2.8573) (1.2075,2.8528) (1.2112,2.8484) (1.2148,2.8440) (1.2184,2.8397) (1.2221,2.8354) (1.2257,2.8311) (1.2293,2.8268) (1.2330,2.8226) (1.2366,2.8184) (1.2402,2.8142) (1.2439,2.8101) (1.2475,2.8059) (1.2511,2.8018) (1.2548,2.7978) (1.2584,2.7937) (1.2620,2.7897) (1.2657,2.7857) (1.2693,2.7818) (1.2729,2.7778) (1.2766,2.7739) (1.2802,2.7700) (1.2838,2.7661) (1.2875,2.7623) (1.2911,2.7585) (1.2947,2.7547) (1.2984,2.7509) (1.3020,2.7472) (1.3056,2.7435) (1.3093,2.7398) (1.3129,2.7361) (1.3165,2.7325) (1.3202,2.7288) (1.3238,2.7252) (1.3274,2.7217) (1.3311,2.7181) (1.3347,2.7146) (1.3383,2.7111) (1.3420,2.7076) (1.3456,2.7041) (1.3492,2.7006) (1.3529,2.6972) (1.3565,2.6938) (1.3602,2.6904) (1.3638,2.6871) (1.3674,2.6837) (1.3711,2.6804) (1.3747,2.6771) (1.3783,2.6738) (1.3820,2.6705) (1.3856,2.6673) (1.3892,2.6641) (1.3929,2.6609) (1.3965,2.6577) (1.4001,2.6545) (1.4038,2.6514) (1.4074,2.6482) (1.4110,2.6451) (1.4147,2.6420) (1.4183,2.6390) (1.4219,2.6359) (1.4256,2.6329) (1.4292,2.6299) (1.4328,2.6269) (1.4365,2.6239) (1.4401,2.6209) (1.4437,2.6180) (1.4474,2.6150) (1.4510,2.6121) (1.4546,2.6092) (1.4583,2.6063) (1.4619,2.6035) (1.4655,2.6006) (1.4692,2.5978) (1.4728,2.5950) (1.4764,2.5922) (1.4801,2.5894) (1.4837,2.5866) (1.4873,2.5839) (1.4910,2.5812) (1.4946,2.5784) (1.4982,2.5757) (1.5019,2.5730) (1.5055,2.5704) (1.5091,2.5677) (1.5128,2.5651) (1.5164,2.5625) (1.5201,2.5598) (1.5237,2.5572) (1.5273,2.5547) (1.5310,2.5521) (1.5346,2.5495) (1.5382,2.5470) (1.5419,2.5445) (1.5455,2.5420) (1.5491,2.5395) (1.5528,2.5370) (1.5564,2.5345) (1.5600,2.5321) (1.5637,2.5296) (1.5673,2.5272) (1.5709,2.5248) (1.5746,2.5224) (1.5782,2.5200) (1.5818,2.5176) (1.5855,2.5153) (1.5891,2.5129) (1.5927,2.5106) (1.5964,2.5083) (1.6000,2.5059)};
\node[anchor=west,font=\small] at (axis cs:1.6,4.810) {$h(e^{-T})$};
\node[anchor=west,font=\small] at (axis cs:1.6,3.870) {$w(T)$};
\node[anchor=west,font=\small] at (axis cs:1.6,2.506) {$v(1-e^{-T})$};
\node[anchor=south,font=\small,gray!70!black] at (axis cs:1.48,1) {$\varkappa(T)+1$};

\end{axis}
\end{tikzpicture}
\caption{
For $T\approx -\log(j/n)$
the relative rank $Y_j$  falling   strictly   below the step curve
 must   be accepted.
The smooth curves are $v(1-e^{-T})<w(T)<h(e^{-T})$.
 }\label{Fig2}
\end{figure}

Thus no stopping occurs before time $(T-\delta_1)_+$,
and the result (\ref{CMRS}) of \cite{CMRS} obtains as
$$w(T)=e^Tu(T)=e^Te^{\delta_1-T}=e^{\delta_1} =3.869\ldots,~~~~T\geq\delta_1.$$
The bounds  on the limit value function
$$
 \frac{2}{1-e^{-T}}=v(1-e^{-T}) <w(T)<h(e^{-T})=e^T \tan\left( \frac{\pi-T}{2}\right)
$$
(where $T<\pi/2$ for the right side; for $T\geq\pi/2$ the upper bound is $h(t_0)=e^{\pi/2}$) are shown in Figure \ref{Fig2}.

{\bf Acknowledgement.} The author is grateful to the two anonymous referees for their careful reading and comments, which led to substantial improvements of the manuscript.


\begin{thebibliography}{99}


\bibitem{SK} Assaf, D. and Samuel-Cahn  E. (1996) The secretary problem: minimizing the expected rank, with i.i.d. random variables,
{\it Adv. Appl. Probab.} {\bf 28}, 828-852.

\bibitem{BDK} Bayraktar, E., Dayanik, S. and Karatzas, I. (2006) Adaptive Poisson disorder problem, {\it Ann. Appl. Probab.} {\bf 16}, 1190--1261.

\bibitem{BG} Berezovsky, B.A. and Gnedin, A.V. {\it The Best Choice Problem}, Moscow, Nauka, 1984.

\bibitem{Brien} O'Brien, G. L. (1983)
Optimal stopping when sampling with and without replacement, {\it Probab. Th. Related Fields}
{\bf 64}, 125--128.

\bibitem{BR} Bruss, F.T. and  Rogers, L.C.G. (1991) Embedding optimal selection problems in a Poisson process,
{\it Stochast. Proc. Appl.} {\bf 38}, 267--278.

\bibitem{BrussSwan} Bruss, F.T. and Swan, Y.C. (2009) A continuous-time approach to Robbins' problem of minimizing the expected rank, {\it J. Appl. Prob.} {\bf 46}, 1--18.

\bibitem{CRS} Chow, Y.S., Robbins, H. and Siegmund, D.  {\it The Theory of Optimal Stopping}, Dover, 1991.

\bibitem{Cayley} Cayley A. (1875) Mathematical questions and their solutions, In: {\it The Collected Mathematical Papers of Arthur Cayley}, vol X (1896), CUP, 587--588.

\bibitem{CMRS} Chow, Y.S., Moriguti, S., Robbins, H. and Samuels, S.M. (1964) {Optimal selection based on relative ranks (the `Secretary Problem')}, {\it Israel J. Math.} {\bf 2}, 81--90.

\bibitem{Ent}
Entwistle H.N.,  Lustri, C.J. and Sofronov G.U. (2022)
On Asymptotics of Optimal Stopping Times,
{\it Mathematics (MDPI)}
{\bf 10}, 194. 

\bibitem{Ferguson} Ferguson, T.  {\it Optimal Stopping and Applications}, e-text, UCLA, 2006.
\bibitem{Finch} Finch, S. (2024) A deceptively simple quadratic recurrence, {\tt arXiv:2409.03510}.

\bibitem{Gianini} Gianini, J. (1977) The infinite secretary problem as the limit of the finite problem, {\it Ann. Probab.} {\bf 5}, 636--644.

\bibitem{GS1} Gianini J.,  and Samuels, S.M. (1976) The infinite secretary problem, {\it Ann. Probab.} {\bf 4}, 418--432.

\bibitem{GnedinSPA} Gnedin, A. (2004) Best choice from the planar Poisson process, {\it Stochastic Proc. Appl.}  {\bf 111}, 317--354.

\bibitem{GC} Gnedin, A. and Carreira, M.C.S. (2026) The memoryless best-choice problem, {\tt arXiv:2607.16145}.

\bibitem{GD1} Gnedin, A. and Derbazi, Z. (2022) Trapping the ultimate success, {\it Mathematics (MDPI) } {\bf 10} (1), 158.

\bibitem{GD2} Gnedin, A. and Derbazi, Z. (2025) The last-success stopping problem with random observation times,
{\it Math. Meth. Operations Research}, {\bf 101}, 1--27.

\bibitem{GJM} Gnedin, A., Janson, S. and Malinovsky, Y. (2026)
Maximal counts in the stopped occupancy problem, 
{\it Electron. J. Probab.} {\bf  31},  1--38. {\tt doi:10.1214/26-EJP1583}


\bibitem{GKS} Gnedin, A., Koziel, P. and Sulkowska M.  (2023)  Running minimum in the best choice problem, {\it Extremes} {\bf 26}, 157--182.

\bibitem{GSe} Gnedin, A. and Seksenbaev, A.  (2021)
Diffusion approximations in the online increasing subsequence problem,
{\it Stoch. Proc.  Appl.} {\bf 139},  298--320.

\bibitem{Malin} Goldenshluger, A., {Malinovsky, Y.} and  Zeevi, A. (2024) Optimal single threshold stopping rules and sharp prophet inequalities,
arXiv:2404.12949

\bibitem{Katriel} Katriel, G. (2025) The Cayley-Moser problem with Poissonian arrival of offers,  
{\tt arXiv:2511.02763}

\bibitem{KK1} Kennedy, D.P. and Kertz, R.P. (1990) Limit theorems for threshold-stopped random variables with applications to optimal stopping, {\it Adv. Appl. Probab.}
{\bf 22}, 396--411.

\bibitem{KK2} Kennedy, D.P. and Kertz, R.P. (1991) The asymptotic behavior of the reward sequence in the optimal stopping of i.i.d. random variables, {\it Ann. Probab.} {\bf 19}, 329--341.

\bibitem{Kolchin} Kolchin, V.F.,
Sevastyanov, B. A. and  Chistyakov, V. P.
{\it  Random Allocations.}
Scripta Series in Mathematics. Halsted Press, New York,  1978.

\bibitem{Lindley} Lindley, D.V. (1961) Dynamic programming and decision theory, {\it Appl. Statistics} {\bf 10}, 39--51.

\bibitem{Mazalov} Mazalov, V.V. and  Peshkov, N.V. (2004) On asymptotic properties of optimal stopping time. {\it Theory Probab.  Appl. } {\bf 48}, 549--555.

\bibitem{Moser} Moser, L. (1956). On a problem of Cayley,
{\it Scripta Math.} {\bf  22}, 289--292.

\bibitem{Mucci1} Mucci, A. G. (1973) {On a class of secretary problems}, {\it Ann. Probab.} {\bf 1}, 417--427.

\bibitem{Mucci2} Mucci, A. G. (1973) { Differential equations and optimal choice problems}, {\it Ann. Statist.} {\bf 1}, 104--113.

\bibitem{Peskir} Peskir, G. and Shiryaev, A. {\it Optimal Stopping and Free-Boundary Problems}, Springer,  2006.

\bibitem{Resnick} Resnick, S. {\it Extreme Values, Regular Variation and Point Processes}, Springer, 2008.

\bibitem{Rusch} R{\"u}schendorf, L. (2016) Approximative solutions of optimal stopping and selection problems,
{\it Mathematica Applicanda} {\bf {22/60}}, 17--44.

\end{thebibliography}
\end{document}